\documentclass[12pt]{article}

\usepackage[english]{babel}

\usepackage{amsfonts, amssymb, amsmath, amsthm, epsf, cite}

\usepackage{graphicx}

\usepackage{fullpage}

\newtheorem{theorem}{Theorem}[section]
\newtheorem{lemma}{Lemma}[section]
\newtheorem{corollary}{Corollary}[section]
\newtheorem{condition}{Condition}[section]

\theoremstyle{definition}
\newtheorem{definition}{Definition}[section]

\newtheorem{example}{Example}[section]
\newtheorem{remark}{Remark}[section]

\numberwithin{equation}{section} \numberwithin{figure}{section}

\renewcommand{\Re}{{\rm Re\,}}

\newcommand{\ph}{{\varphi}}

\newcommand{\cA}{{\mathcal A}}
\newcommand{\cD}{{\mathcal D}}

\newcommand{\cT}{{\mathcal T}}

\newcommand{\cP}{{\mathcal P}}

\newcommand{\cB}{{\mathcal B}}

\newcommand{\cR}{{\mathcal R}}

\newcommand{\cU}{{\mathcal U}}

\newcommand{\ccQ}{{\mathcal Q}}

\newcommand{\bdelta}{{\boldsymbol\delta}}
\newcommand{\bDelta}{{\boldsymbol\Delta}}

\newcommand{\bA}{{\mathbf A}}

\newcommand{\bF}{{\mathbf F}}

\newcommand{\bu}{{\mathbf u}}
\newcommand{\br}{{\mathbf r}}
\newcommand{\bv}{{\mathbf v}}

\newcommand{\bbR}{{\mathbb R}}
\newcommand{\bbN}{{\mathbb N}}
\newcommand{\bbC}{{\mathbb C}}
\newcommand{\bbZ}{{\mathbb Z}}

\newcommand{\ep}{\varepsilon}

\title{Asymptotics of parabolic Green's functions on lattices}

\author{Pavel  Gurevich\footnote{Free University of Berlin, Germany; Peoples' Friendship
University of Russia, Russia; email: gurevich@math.fu-berlin.de}}

\begin{document}

\maketitle

\begin{abstract}
For parabolic spatially discrete equations, we consider Green's functions, also
known as heat kernels on lattices. We obtain their asymptotic
expansions with respect to powers of time variable $t$ up to
an arbitrary order and estimate the remainders uniformly on the
whole lattice. The spatially discrete (difference) operators under consideration are
finite-difference approximations of continuous strongly elliptic
differential operators (with constant coefficients) of arbitrary
even order in $\bbR^d$ with arbitrary $d\in\bbN$. This genericity,
besides numerical and deterministic lattice-dynamics applications,
allows one to obtain higher-order asymptotics of transition
probability functions for continuous-time random walks on $\bbZ^d$ and other lattices.
\end{abstract}



\tableofcontents

\section{Introduction}

The paper deals with Green's functions, or heat kernels, of
general parabolic equations with constant coefficients that are
continuous in time and discrete in space. Applications that we
have in mind include spatial discretizations of continuous models, lattice dynamical systems, and continuous-time random walks.

We consider parabolic problems on the {\it grid space}, or {\it
lattice},
$$
\bbR_\ep^d:=\{x\in\bbR^d: x_k=s_k\ep,\ s_k\in\bbZ,\
k=1,\dots,d\},\quad d\in\bbN,\ \ep>0,
$$
of the form
\begin{equation}\label{eqz_nGreenMultiDimIntro}
\left\{
\begin{aligned}
& \dot \bu^\ep(x,t)+\bA_\ep \bu^\ep(x,t)=0, & & x\in\bbR_\ep^d,\ t>0,\\
& \bu^\ep(x,0)=\bdelta^\ep(x),   & & x\in\bbR_\ep^d.
\end{aligned}
\right.
\end{equation}
Here $\dot{}=\partial/\partial t$, $\bdelta^\ep(x)$ is the {\it
grid delta-function} given by
\begin{equation}\label{eqdeltaEp}
\bdelta^\ep(0)=\ep^{-d},\qquad \bdelta^\ep(x)=0\quad \forall
x\in\bbR_\ep^d\setminus\{0\}.
\end{equation}
and $\bA_\ep$ is an elliptic difference operator (with constant
coefficients) which is assumed to be an $M$th order approximation
of a strongly elliptic (continuous) differential operator
$\cA(\cD)$ of even order $\ell$ (see
Section~\ref{secHigherDAsymptotics} for rigorous definitions). An
order $M\in\bbN$ is defined by the estimate
\begin{equation}\label{eqMApproxIntro}
\sup\limits_{y\in\bbR^d}|\bA_\ep u(y)-\cA(\cD)u(y)|\le
C(u)\ep^M\quad \forall \ep>0,
\end{equation}
which must hold for any {\it smooth} (continuous and bounded with all
its derivatives) function $u(y)$ and an appropriate constant
$C(u)>0$ not depending on $\ep>0$.

We call a solution of problem~\eqref{eqz_nGreenMultiDimIntro} the
{\it
 first discrete Green function}. Given  arbitrary initial data
instead of $\bdelta^\ep(x)$, one can represent a solution as a
(discrete) convolution of these initial data with the first Green
function. Alternatively, the first Green function determines a
semigroup generated by $-\bA_\ep$. Various properties of $\bA_\ep$
(in terms of its symbol) in functional spaces of grid functions,
corresponding semigroup properties, and relations to spatially
continuous problems can be found in~\cite{Frank,Sobolevskii,Beyn}, see
also references therein.

By using the discrete Fourier transform, one can give an integral
representation of the discrete Green function (see~\eqref{eqz_nFourierMultiDim}). However, it is not
possible to express it via elementary functions, hence its
asymptotic expansions play an important role. Asymptotics of the
so-called lattice
 Green's functions in the stationary case were studied beginning from 1950s, see~\cite{Duffin}, the subsequent papers~\cite{Bramble, Katsura,
Mangad, Martinsson, Guttmann, Otto}, and the monograph~\cite[Chapter~8]{Lifanov}.
For parabolic operators, there is vast literature in the spatially continuous case. For example, large-time behavior of Green's functions was treated in~\cite{Murata85} (for small perturbations of the heat operator) and in \cite{Tsuchida08,Norris97} (for spatially periodic coefficients). A survey on the large time behavior of heat kernels for second-order parabolic operators on Riemannian manifolds can be found in~\cite{Pinchover13}. In the spatially discrete case, the research directions include
continuous-time random walks on general graphs (see, e.g.,~\cite{Lawler,
Pang93} and references therein) and on lattices in a random
environment (see, e.g.,~\cite{Delmotte}). In both cases,
Gaussian bounds for the heat kernel is an important question. However, higher-order asymptotics of  Green's functions is not available in general. We mention~\cite{Iliev02, Iliev07} , where an asymptotic expansion of  Green's function for specific parabolic equations on one-dimensional lattices was obtained in terms of the Bessel functions.

In our paper, using the integral Fourier representation of the
first Green function $\bu^\ep(x,t)$, we obtain a higher-order
asymptotic formula of the form
\begin{equation}\label{eqAsymphGeneralMultiDimIntro}
\begin{aligned}
\dfrac{\partial^J\bu^\ep(x,t)}{\partial t^J} &
=\dfrac{1}{t^{d/\ell+J}}\,
{H}_{J}\left(\dfrac{x}{t^{1/\ell}}\right)\\
& + \sum\limits_{k=M}^K \dfrac{\ep^k}{t^{(k+d)/\ell+J}}\,
{H}_{Jk}\left(\dfrac{x}{t^{1/\ell}}\right) + \br_\bu^\ep(J,K;
x,t),\quad x\in\bbR_\ep^d,\ t\ge t_0\ep^\ell.
\end{aligned}
\end{equation}
Here $\ep>0$, $t_0>0$, $K\ge M$, $J\ge 0$, ${H}_{J}(y)$ and
${H}_{Jk}(y)$, $y\in\bbR^d$, are explicitly given smooth
functions, and the remainder satisfies
\begin{equation}\label{eqHkRemainderIntro}
|\br_\bu^\ep(J,K; x,t)|\le
\dfrac{\ep^{K+1}R_\bu(J,K,t_0)}{t^{(K+d+1)/\ell+J}}, \quad
x\in\bbR_\ep^d,\ t\ge t_0\ep^\ell,
\end{equation}
with $R_\bu(J,K,t_0)\ge 0$ not depending on $\ep>0$,
$x\in\bbR_\ep^d$, and $t\ge t_0\ep^\ell$. We emphasize that the
estimate in~\eqref{eqHkRemainderIntro} is uniform with respect to
$x\in\bbR_\ep^d$.

Moreover, we show that, for $J=0$, the leading order term
$\dfrac{1}{t^{d/\ell}}\,
{H}_{0}\left(\dfrac{x}{t^{1/\ell}}\right)$
in~\eqref{eqAsymphGeneralMultiDimIntro} coincides with the
continuous parabolic Green function (see~\eqref{eqContGreen}).
Hence, \eqref{eqAsymphGeneralMultiDimIntro} implies that the first
discrete Green function approximates the continuous one, uniformly for $x\in\bbR_\ep^d$, with an
error $O\left(\dfrac{1}{t^{(M+d)/\ell}}\right)$ as $t\to\infty$
and with an error $O(\ep^M)$ as $\ep\to0$, where $M$ is the order
of approximation in~\eqref{eqMApproxIntro} (see
Corollaries~\ref{corApproxt} and~\ref{corApproxEpsilon}).
 In the
case where $-\bA_\ep$ is an approximation of the Laplacian, the
first term in~\eqref{eqAsymphGeneralMultiDimIntro} with $J=0$
equals $ \dfrac{1}{(2\sqrt{\pi t})^d} e^{-|x|^2/4t}$. In particular, it is radially symmetric.
The other terms, being non-symmetric,
take into account the radial non-symmetry of the lattice
$\bbR_\ep^d$.

Asymptotic formula~\eqref{eqAsymphGeneralMultiDimIntro} is proved
in Section~\ref{secFirstGreen} (Theorem~\ref{thAsymphMultiDim}), after rigorous definitions of elliptic differential and difference operators in
Section~\ref{secHigherDAsymptotics}.

We note that the general form of the difference operator $\bA_\ep$
(rather than the Laplacian only) allows one to treat
continuous-time random walks on $d$-dimensional lattices that are not necessarily cubic, but still represented by a discrete additive subgroup of $\bbR^d$. This is
possible whenever one can linearly transform the vertices of such a lattice
to $\bbR_\ep^d$ (see~\cite[Section~1]{Lawler}) and obtain an
elliptic generator~$\bA_\ep$ on $\bbR_\ep^d$. Then,
linearly transforming~\eqref{eqAsymphGeneralMultiDimIntro} back to the
original non-cubic lattice, one can deduce asymptotics of the transition
probability function of the original random walk.
Example~\ref{exTwoLattices}.2 illustrates an elliptic difference
operator that is a generator obtained after transforming the
random walk on the triangular lattice to $\bbZ^2$.

Along with the first discrete Green function, we study the {\it
second discrete Green function} in
Sections~\ref{secSecondGreenPrelim} and~\ref{secSecondGreen}. It
is defined as a solution of the problem
\begin{equation}\label{eqy_nGreenMultiDimIntro}
\left\{
\begin{aligned}
& \dot \bv^\ep(x,t)+\bA_\ep \bv(x,t)=\bdelta^\ep(x), & & x\in\bbR_\ep^d,\ t>0,\\
& \bv^\ep(x,0)=0,   & & x\in\bbR_\ep^d.
\end{aligned}
\right.
\end{equation}
Given an arbitrary time-independent right-hand side instead of
$\bdelta^\ep(x)$, one can represent a solution as a (discrete)
convolution of this right-hand side with the second Green
function. An interesting application of the second Green function $\bv^\ep(x,t)$ occurs in discrete reaction-diffusion equations with hysteretic nonlinearity in the right-hand side. Higher-order asymptotic formulas for $\bv^\ep(x,t)$ play a central role in understanding pattern formation mechanisms there (see~\cite{GurTikhRattling2015}).

Asymptotics of the second Green function depends on the sign of
$\ell-d$, where $\ell$ is the order of the continuous differential
operator $\cA(\cD)$ and $d$ is the spatial dimension. If
$\ell-d\le -1$, then we prove that  the second discrete Green
function satisfies
\begin{equation}\label{eqAsympfGeneralMultiDim-d-lBigger1Intro}
\begin{aligned}
\bv^\ep(x,t)&=  \dfrac{1}{t^{d/\ell-1}} \,
{F}_0\left(\dfrac{x}{t^{1/\ell}}\right) + \dfrac{1}{\ep^{d-\ell}}\Omega\left(\dfrac{x}{\ep}\right)\\
&+\sum\limits_{k=M}^K\dfrac{\ep^k}{t^{(k+d)/\ell-1}}\,
{F}_k\left(\dfrac{x}{t^{1/\ell}}\right) + \br_\bv^\ep(K; x,t),
\quad x\in\bbR_\ep^d,\ t\ge t_0\ep^\ell,
\end{aligned}
\end{equation}
where the remainder is estimated as
\begin{equation}\label{eqRemainderIntro}
|\br_\bv^\ep(K; x,t)|\le
\dfrac{\ep^{K+1}R_\bv(K,t_0)}{t^{(K+d+1)/\ell-1}}, \quad
x\in\bbR_\ep^d,\ t\ge t_0\ep^\ell,
\end{equation}
with $R_\bv(K,t_0)\ge 0$ not depending on $\ep>0$,
$x\in\bbR_\ep^d$, and $t\ge t_0\ep^\ell$. Again, the
estimate in~\eqref{eqHkRemainderIntro} is uniform with respect to
$x\in\bbR_\ep^d$. In
Section~\ref{secSecondGreenPrelim}
(Lemma~\ref{lAsymphMultiDimPrelim}), we
prove~\eqref{eqAsympfGeneralMultiDim-d-lBigger1Intro} with the
functions ${F}_k(y)$ that are given explicitly  as solutions of
certain first-order PDEs with ${H}_0(y)$ and ${H}_{0k}(y)$ in the
right-hand side, but still unknown function $\Omega(y)$. In
Section~\ref{subsecEll-dNegative}
(Theorem~\ref{thAsympfGeneral1}), we find the function $\Omega(y)$
by simultaneously passing to the limit, as $t\to\infty$,
in~\eqref{eqAsympfGeneralMultiDim-d-lBigger1Intro} and in the
explicit integral representation of $\bv^\ep(x,t)$.

If $\ell-d\ge 0$, the situation is more delicate. In
Section~\ref{secSecondGreenPrelim}, we prove an analogue
of~\eqref{eqAsympfGeneralMultiDim-d-lBigger1Intro} in the form
\begin{equation}\label{eqAsympfGeneralMultiDim-d-lBigger1Final1Intro}
\begin{aligned}
\bv^\ep(x,t) & = t^{1-d/\ell}\,
{F}_0\left(\dfrac{x}{t^{1/\ell}}\right) +
\sum\limits_{k=M}^{\ell-d} \ep^k t^{1-(k+d)/\ell}\,
{F}_k\left(\dfrac{x}{t^{1/\ell}}\right) + \ep^{\ell-d}\,
\Omega\left(\dfrac{x}{\ep}\right) \\
& + \sum\limits_{k=\max(M,\ell-d+1)}^K \dfrac{\ep^k}{
t^{(k+d)/\ell-1}}\, {F}_k\left(\dfrac{x}{t^{1/\ell}}\right) +
\br_\bv^\ep(K; x,t), \quad x\in\bbR_\ep^d\setminus\{0\},\ t\ge
t_0\ep^\ell
\end{aligned}
\end{equation}
(if $M\ge l-d+1$, the first sum is absent). As before,  the remainder
satisfies~\eqref{eqRemainderIntro}, the functions ${F}_k(y)$
  are given explicitly  as solutions of certain first-order
PDEs with ${H}_0(y)$ and ${H}_{0k}(y)$ in the right-hand side, and
$\Omega(y)$ is still unknown. However, unlike in the case
$\ell-d\le -1$,
formula~\eqref{eqAsympfGeneralMultiDim-d-lBigger1Final1Intro} may
contain the function ${F}_{\ell-d}(y)$, which in general appears to be
undefined at $y=0$ and to have logarithmic growth as $y\to 0$.
Furthermore, we cannot pass to the limit as $t\to\infty$
immediately because both the terms with $k=0,M,\dots,\ell-d$
in~\eqref{eqAsympfGeneralMultiDim-d-lBigger1Final1Intro} and the
explicit integral representation of $\bv^\ep(x,t)$ tend to
infinity. Hence, in Section~\ref{subsecEll-dPositive}
(Theorem~\ref{thAsympfGeneral2}), we use the explicit integral
representation of $\bv^\ep(x,t)$ to deduce another asymptotic
representation. It contains a linear combination of nonnegative
powers $t^{1-(k+d)/\ell}$ ($k=0,\dots,\ell-d$) and of $\ln t$ with
$x$-dependent coefficients and holds for each fixed
$x\in\bbR_\ep^d$ (but not uniformly in $\bbR_\ep^d$). Comparing it
with~\eqref{eqAsympfGeneralMultiDim-d-lBigger1Final1Intro} allows
one to determine $\Omega(y)$ and obtain uniform asymptotics in~$\bbR_\ep^d$.

In Section~\ref{sec1DExample}, we apply our general results to the
$(2N)$th order approximation of the one-dimensional Laplacian:
\begin{equation}\label{eqDeltaEpIntro}
\bA_\ep=-\bDelta_\ep,\quad \bDelta_\ep
\bu(x):=\ep^{-2}\sum\limits_{\nu=1}^N
a_\nu\big(\bu(x-\ep\nu)-2\bu(x)+\bu(x+\ep\nu)\big)
\end{equation}
with appropriately chosen coefficients $a_\nu\in\bbR$. (In this
case, $\ell=2$ and $d=1$, hence we are in the situation $\ell-d\ge
0$.) First, we claim that $\bA_\ep$ is elliptic for any $N\in\bbN$
in the sense of Condition~\ref{condDiscreteEllipticity} (the proof
is  given in Appendix~\ref{appendProofEllipt}). Then we
explicitly find all the functions in
expansions~\eqref{eqAsymphGeneralMultiDimIntro}
and~\eqref{eqAsympfGeneralMultiDim-d-lBigger1Final1Intro} for the
first and second Green functions and in the corresponding
expansions for the (spatial) gradients of the first and second
Green functions.

Interestingly, if $N=1$, i.e., $\bDelta_\ep
\bu(x):=\ep^{-2}\big(\bu(x-\ep)-2\bu(x)+\bu(x+\ep)\big)$, it appears that
$\Omega(x/\ep)$
in~\eqref{eqAsympfGeneralMultiDim-d-lBigger1Final1Intro} vanishes
for all $x\in\bbR_\ep^d$, which is proved in
Appendix~\ref{appendGamma0}. In general, $\Omega(x/\ep)\not\equiv
0$ for $N\ge 2$.

To conclude, we collect some general notation that we use throughout the paper. We use variables $x\in\bbR_\ep^d$ or
$x\in\bbZ^d$ (for $\ep=1$) and $y\in\bbR^d$ for the grid and
continuous coordinates, respectively, and $\eta\in
R_{\pi\ep^{-1}}$, $\theta\in  R_{\pi}$, and $\xi\in\bbR^d$ for the
Fourier coordinates, where
\begin{equation}\label{eqRepsilon}
R_{\lambda}:=\{\xi\in\bbR^d: |\xi_k|\le \lambda,\
k=1,\dots,d\},\quad  \lambda>0.
\end{equation}
By $B_\lambda$, we denote the ball of radius $\lambda>0$ in $\bbR^d$.
We denote by $(r,\varphi_1,\dots\varphi_{d-1})$ the spherical
coordinates of $y\in\bbR^d$ or $\theta\in R_\pi$ and introduce the function
\begin{equation}\label{eqr_pi}
r_\pi(\ph):=\text{distance from the origin to the boundary of
 $R_\pi$ in the direction $\ph$}.
\end{equation}
We use the upper-case bold letters
 for  difference operators: $\bA_\ep$, $
\bDelta_\ep$ and lower-case bold letters
for grid functions (i.e., functions defined on~$\bbR_\ep^d$):
$\bu^\ep(x,t)$, $\bv^\ep(x,t)$, etc., indicating by superscript
$\ep$ that $x\in\bbR_\ep^d$. We may omit the superscript if
$\ep=1$ and $x\in\bbZ^d$. We denote polynomials by calligraphic letters:
$\cA(\xi)$, $\cP(\xi)$, etc. and other (non-polynomial) functions by usual letters: $A(\theta)$, $H(y)$, $\Omega(y)$,
$u(y,t)$, etc.

\section{Continuous and discrete parabolic Green's functions}\label{secHigherDAsymptotics}

In this section, we introduce parabolic equations with general
 elliptic operators of even order~$\ell$ in $\bbR^d$, $d\in\bbN$, and the corresponding finite-difference approximation and define
their Green functions.

\subsection{Continuous Green function}
Consider a differential operator with constant coefficients in
$\bbR^d$, $d\in\bbN$:
\begin{equation}\label{eqcA}
\cA(\cD)=\sum\limits_{|\alpha|=\ell}b_\alpha \cD^\alpha,
\end{equation}
where $\ell\in\bbN$ is even, $\alpha=(\alpha_1,\dots,\alpha_d)$ is
a multi-index, $|\alpha|=\alpha_1+\dots+\alpha_d$,
$b_\alpha\in\bbC$,
$\cD^\alpha=\cD_1^{\alpha_1}\dots\cD_d^{\alpha_d}$, and
$\cD_j=-i\partial/\partial y_j$, $j=1,\dots,d$.

\begin{definition}
The polynomial
\begin{equation}\label{eqcAtheta}
\cA(\xi):=\sum\limits_{|\alpha|=\ell}b_\alpha \xi^\alpha,\quad
\xi\in\bbR^d,
\end{equation}
where $\xi^\alpha=\xi_1^{\alpha_1}\cdot\dots\cdot
\xi_d^{\alpha_d}$ (obtained from $\cA(\cD)$ by formally replacing
$\cD_j$ by $\xi_j$), is called the {\em symbol} of the differential operator
$\cA(\cD)$.
\end{definition}

Obviously, there is $C>0$ such that
\begin{equation}\label{eqcAUpperBound}
|\cA(\xi)|\le C|\xi|^\ell\quad \forall\xi\in\bbR^d.
\end{equation}

We assume throughout that the operator $\cA(\cD)$ is {\em strongly
elliptic}, i.e., the following condition holds.
\begin{condition}\label{condAContElliptic}
There is $c>0$ such that
\begin{equation}\label{eqcAelliptic}
\Re\cA(\xi)\ge c|\xi|^\ell\quad \forall\xi\in\bbR^d.
\end{equation}
\end{condition}

We briefly recall the notion of a {\it continuous Green function}
for   parabolic equations. We define it as a solution of the
problem
\begin{equation}\label{eqParabolcA}
\left\{
\begin{aligned}
&\dfrac{\partial u(y,t)}{\partial t}+\cA(\cD) u(y,t)=0, & & y\in \bbR^d,\ t>0,\\
&u(y,0)=\delta(y),& & y\in \bbR^d,
\end{aligned}
\right.
\end{equation}
with the $\delta$-function in the initial condition. To give an
explicit formula for the Green function, we formally use the
Fourier transform and deduce from~\eqref{eqParabolcA}
$$
u(y,t)=\dfrac{1}{(2\pi)^d}\int_{\bbR^d} e^{-t\cA(\xi)}e^{i
y\xi}\,d\xi,\quad y\in\bbR^d,\ t>0,
$$
where the integral converges for each fixed $t>0$ due
to~\eqref{eqcAelliptic}. Making the change of variables
$\xi\mapsto\xi t^{-1/\ell}$ and using the  homogeneity of $\cA(\xi)$, we obtain
\begin{equation}\label{eqContGreen}
u(y,t)=\dfrac{1}{t^{d/\ell}}
{H}\left(\dfrac{y}{t^{1/\ell}}\right),\quad y\in\bbR^d,\ t>0,
\end{equation}
where
\begin{equation}\label{eqhGeneral}
{H}(y):= \dfrac{1}{(2\pi)^d}\int_{\bbR^d} e^{-\cA(\xi)}e^{i
y\xi}\,d\xi.
\end{equation}
Due to~\eqref{eqcAelliptic}, ${H}(y)$ and all its derivatives
decay at infinity faster than any negative power of~$|y|$.

\begin{definition}
We call the function given by~\eqref{eqContGreen} the {\em
continuous Green function}.
\end{definition}

\begin{remark} If $\cA(\cD)=-\Delta$, where $\Delta$ is the
Laplace operator, then $\cA(\xi)=|\xi|^2$, $ {H}(y)=\dfrac{1}{2^d\pi^{d/2}}
e^{-|y|^2/4} $, and~\eqref{eqContGreen} assumes the well-known form
$u(y,t)=\dfrac{1}{(2\sqrt{\pi t})^d} e^{-|y|^2/4t}$.
\end{remark}

\subsection{Discrete Green function}
For $\ep>0$, we
define the {\it grid space}
$$
\bbR_\ep^d:=\{x\in\bbR^d: x_k=s_k\ep,\ s_k\in\bbZ,\ k=1,\dots,d\}.
$$
A function defined on $\bbR_\ep^d$ is called a {\it grid
function}. We say that a grid function (continuous function) is
{\em rapidly decreasing} if it decays at infinity faster than any
negative power of $|x|$, $x\in\bbR_\ep^d$ (of $|y|$,
$y\in\bbR^d$).
For a continuous function $u(y)$, $y\in\bbR^d$, we use the
notation
\begin{equation}\label{eqFiniteDifference}
\delta_{\ep,k\pm}u(y)=u(y\pm e_k\ep)-u(y),
\end{equation}
where $e_k$ is the unit vector of the axis $y_k$. The same notation is used for grid functions $\bu^\ep(x)$, $x\in\bbR_\ep^d$.

We fix a natural $\hat \ell\ge \ell$ and consider the following
difference operator with constant coefficients:
\begin{equation}\label{eqbA}
\bA_\ep:=\ep^{-\ell}\sum\limits_{\ell\le|\nu|\le \hat \ell} a_\nu
\delta_{\ep,1-}^{\nu_1}\delta_{\ep,1+}^{\nu_2}\dots
\delta_{\ep,d-}^{\nu_{2d-1}}\delta_{\ep,d+}^{\nu_{2d}},
\end{equation}
where $\nu=(\nu_1,\dots,\nu_{2d})$ is a multi-index,
$|\nu|=\nu_1+\dots+\nu_{2d}$, and $a_\nu\in\bbC$.

We say that a function $u(y)$, $y\in\bbR^d$, is {\it smooth} if it
is continuous and bounded together with all its derivatives.

\begin{definition}\label{defAppriximation}
For $M\in\bbN$, the difference operator $\bA_\ep$ is said to be an
{\em $M$th order approximation of the differential operator
$\cA(\cD)$} if, for any smooth function $u(y)$, there is a
constant $C(u)>0$ such that
\begin{equation}\label{eqMApprox}
\sup\limits_{y\in\bbR^d}|\bA_\ep u(y)-\cA(\cD)u(y)|\le
C(u)\ep^M\quad \forall \ep>0.
\end{equation}
\end{definition}

We assume throughout that the following holds.
\begin{condition}\label{condAppriximation}
There is $M\in\bbN$ such that $\bA_\ep$ is   an   $M$th order
approximation of the differential operator $\cA(\cD)$.
\end{condition}

Next, we introduce a symbol of $\bA_\ep$ and relate it to the
symbol of $\cA(\cD)$. To motivate the notion of a symbol,  we define the Fourier transform of a rapidly decreasing grid
function $\bu^\ep(x)$ by the formula
\begin{equation}\label{eqFourier}
(\bF_\ep \bu^\ep)(\eta):=(2\pi)^{-d}\sum\limits_{x\in\bbR_\ep^d}e^{-i
x\eta}\bu(x)\ep^d, \quad \eta\in R_{\pi\ep^{-1}},
\end{equation}
where $R_{\pi\ep^{-1}}$ is given by~\eqref{eqRepsilon}. The
Fourier transform $\bF_\ep$ establishes an isomorphism between
rapidly decreasing grid functions and smooth functions that are
$2\pi \ep^{-1}$-periodic with respect to each of the variables
$\eta_1,\dots,\eta_d$. The inverse transform is given by
\begin{equation}\label{eqInverseFourier}
(\bF_\ep^{-1} v)(x)=\int_{R_{\pi\ep^{-1}}} e^{i
x\eta}v(\eta)\,d\eta.
\end{equation}

If $\bu^\ep(x)$ is a rapidly decreasing grid function, then $\bA_\ep
\bu^\ep(x)$ is also rapidly decreasing and
\begin{equation}\label{eqFourierSymbol}
(\bF_\ep \bA_\ep \bu^\ep)(\eta) = \ep^{-\ell}A(\eta\ep)\cdot
(\bF_\ep \bu^\ep)(\eta),
\end{equation}
where the function $\ep^{-\ell}A(\eta\ep)$ is obtained by
replacing the operator $\delta_{\ep,k\pm}$ in the right-hand side
of equality~\eqref{eqbA} by the expression $e^{\pm i\eta_k\ep}-1$.

\begin{definition}\label{defSymbolDiscrete}
The function
\begin{equation}\label{eqSymbolDiscrete}
A(\theta):= \sum\limits_{\ell\le|\nu|\le \hat \ell} a_\nu
\left(e^{-i\theta_1}-1\right)^{\nu_1} \left(e^{i\theta_1}-1\right)^{\nu_2}\dots
\left(e^{-i\theta_d}-1\right)^{\nu_{2d-1}} \left(e^{i\theta_d}-1\right)^{\nu_{2d}},\quad \theta\in R_\pi,
\end{equation}
is called the {\em symbol of the difference operator~$\bA_\ep$.}
\end{definition}

The following lemma establishes a connection between the symbols
of $\cA(\cD)$ and $\bA_\ep$.

\begin{lemma}\label{lConnectionSymbols}
If Condition~$\ref{condAppriximation}$ holds, then, for any $K\ge
M$,
\begin{equation}\label{eqConnectionSymbols}
A(\theta)=\sum\limits_{k=0}^K
\cA_{k+\ell}(\theta)+O(|\theta|^{K+\ell+1}),
\end{equation}
where
\begin{equation}\label{eqAk}
\cA_{k+\ell}(\theta):=\left\{
\begin{aligned}
&\cA(\theta) & &\text{if } k=0,\\
&0   & &\text{if } k=1,\dots,M-1,\\
&\text{homogeneous polynomial}\\
&\text{of degree $k+\ell$} & &\text{if } k=M,\dots,K,
\end{aligned}
\right.
\end{equation}
and $O(\cdot)$ is understood with respect to $|\theta|\to0$.

 If the symbol of $\bA_\ep$ satisfies~\eqref{eqConnectionSymbols} and~\eqref{eqAk} with some $K\ge M$, then
 Condition~$\ref{condAppriximation}$ holds.
\end{lemma}
\proof We observe that, for any smooth $u(y)$ and any $J\ge 1$,
$$
\begin{aligned}
\delta_{\ep,k\pm} u(y)&=\sum\limits_{j=1}^J \dfrac{(\pm 1)^j \ep^j}{j!}\dfrac{\partial^j u(y)}{\partial y_k^j} + \ep^{J+1} R_{k,J+1}(y),\\
e^{\pm i\theta_k}-1&=\sum\limits_{j=1}^J \dfrac{(\pm 1)^j}{j!}(i\theta_k)^j + O(|\theta|^{J+1}),\\
\end{aligned}
$$
where $\sup\limits_{y\in\bbR^d}|R_{k,J+1}(y)|\le C_1(u)$ with
some $C_1(u)\ge 0$. Therefore, taking into account definitions~\eqref{eqbA} and~\eqref{eqSymbolDiscrete}, we have
for any $K\ge M$
$$
\begin{aligned}
\bA_\ep u(y)&=\sum\limits_{k=0}^K \ep^k \cA_{k+\ell}(\cD)u(y) + \ep^{K+1} R_{K+\ell+1}(y),\\
A(\theta)&=\sum\limits_{k=0}^K\cA_{k+\ell}(\theta) + O(|\theta|^{K+\ell+1}),\\
\end{aligned}
$$
where $\cA_{k+\ell}(\cdot)$ are the same homogeneous polynomials of degree
$k+\ell$ in both formulas and $\sup\limits_{y\in\bbR^d}|R_{K+\ell+1}(y)|\le
C_2(u)$ with some $C_2(u)\ge 0$. Comparing the last two equalities
and taking into account~\eqref{eqMApprox}, we conclude the proof.
\endproof

Along with Condition~\ref{condAppriximation}, we assume throughout
that the following {\em ellipticity condition} holds.
\begin{condition}\label{condDiscreteEllipticity}
$ \Re A(\theta)>0$ for all $\theta\in R_\pi\setminus\{0\}$, where
$R_\pi$ is defined in~$\eqref{eqRepsilon}$.
\end{condition}

\begin{remark}
Condition~\ref{condDiscreteEllipticity} does not automatically
follow from the fact that $\bA_\ep$ is an approximation of a
strongly elliptic operator $\cA(\cD)$ (only the inequality $\Re A(\theta)>0$ for small $|\theta|\ne0$ follows). However,
Lemma~\ref{lAtheta} below implies that
Condition~\ref{condDiscreteEllipticity}  holds, for example, for
any order approximation of the second derivative in~$\bbR$
(see~\eqref{eqDeltaEp} and~\eqref{eqanu}) and, hence, for the
corresponding approximation of the Laplace operator in $\bbR^d$, $d\in\bbN$.
\end{remark}

\begin{example}\label{exTwoLattices} 1.  Consider the difference operator
$$
\bA_\ep:=\ep^{-2}\sum\limits_{j=1}^d\delta_{\ep,j+}\delta_{\ep,j-}.
$$
For $\ep=1$, the operator $-(2d)^{-1}\bA_1$ is a generator of the continuous-time
random walk on $\bbZ^d$ (see~\cite[Section~1]{Lawler}). The
operator $\bA_\ep$ is the 2nd order approximation of the strongly
elliptic operator $\cA(\cD)=-\Delta$, where $\Delta$
is the Laplacian. The symbol of $\bA_\ep$  is given
by
$$
A(\theta):=2\sum\limits_{j=1}^d(1-\cos \theta_j).
$$
In Section~\ref{sec1DExample}, we will also consider $(2N)$th
order approximations of the Laplacian for $d=1$.

2. Let $d=2$. Consider the operator
\begin{equation}\label{eqLaplaceTriangular}
\bA_\ep:=\dfrac{2}{3\ep^2}\big(2\delta_{\ep,1+}\delta_{\ep,1-}+2\delta_{\ep,2+}\delta_{\ep,2-}-
\delta_{\ep,1+}\delta_{\ep,2-}-
\delta_{\ep,2+}\delta_{\ep,1-}\big).
\end{equation}
The operator $\bA_\ep$
 is the 2nd order approximation of the strongly elliptic operator $\cA(\cD)u:= -\dfrac{4}{3}(u_{y_1y_1}+u_{y_2y_2}-u_{y_1y_2})$. The symbol of $\bA_\ep$ is given by
$$
A(\theta)=4-\dfrac{4}{3}\big(\cos\theta_1+\cos\theta_2+\cos(\theta_1-\theta_2)\big)
$$
and the symbol of $\cA(\cD)$ is given by
$$\cA(\theta)=\dfrac{2}{3}(\theta_1^2+\theta_2^2+(\theta_1-\theta_2)^2).$$

On the other hand, $-\bA_\ep$ corresponds to the discretization on the triangular lattice of the continuous Laplacian. Indeed, the matrix
$$
M:=\begin{pmatrix}
  1 & -1/\sqrt{3} \\
  0 & 2/\sqrt{3}
\end{pmatrix}.
$$
takes the triangular lattice to $\bbZ^2$ and, at the same time, $\cA((M^T)^{-1}\xi)=\xi_1^2+\xi_2^2$. The latter is the symbol of the differential operator $-\Delta$. 
\end{example}

Lemma~\ref{lConnectionSymbols}, together with  Conditions~$\ref{condAContElliptic}$,
$\ref{condAppriximation}$,
and~$\ref{condDiscreteEllipticity}$, implies the following.

\begin{lemma} We have
\begin{equation}\label{eqDiscreteEllipticitymUpperBound}
|A(\theta)|\le C|\theta|^\ell\quad \forall \theta\in R_\pi,
\end{equation}
\begin{equation}\label{eqDiscreteEllipticitym}
\Re A(\theta)\ge c|\theta|^\ell\quad \forall \theta\in R_\pi,
\end{equation}
where $C$ and $c$ are positive constants that do not depend on
$\theta\in R_\pi$.
\end{lemma}

Without loss of generality, we assume that the constants $C$ and $c$ in~\eqref{eqDiscreteEllipticitymUpperBound}
and~\eqref{eqDiscreteEllipticitym}
are the same as in~\eqref{eqcAUpperBound}
and~\eqref{eqcAelliptic}, respectively.

Now we   introduce discrete Green functions. Let $\bdelta^\ep(x)$
be the {\it grid delta-function}, see~\eqref{eqdeltaEp}.

\begin{definition}\label{defGreenMultiDim}
We call the function $\bu^\ep(x,t)$, $x\in\bbR_\ep^d,\ t\ge 0$,
the {\em first discrete Green function} if $\bu^\ep(\cdot,t)$ is a
rapidly decreasing grid function for all $t\ge0$,
$\bu^\ep(x,\cdot)\in C^1[0,\infty)$ for all $x\in\bbR_\ep^d$, and
\begin{equation}\label{eqz_nGreenMultiDim}
\left\{
\begin{aligned}
& \dot \bu^\ep(x,t)+\bA_\ep \bu^\ep(x,t)=0, & & x\in\bbR_\ep^d,\ t>0,\\
& \bu^\ep(x,0)=\bdelta^\ep(x),   & & x\in\bbR_\ep^d.
\end{aligned}
\right.
\end{equation}

We call the function $\bv^\ep(x,t)$, $x\in\bbR_\ep^d,\ t\ge 0$,
the {\em second discrete Green function} if $\bv^\ep(\cdot,t)$ is
a rapidly decreasing grid function for all $t\ge0$,
$\bv^\ep(x,\cdot)\in C^1[0,\infty)$ for all $x\in\bbR_\ep^d$, and
\begin{equation}\label{eqy_nGreenMultiDim}
\left\{
\begin{aligned}
& \dot \bv^\ep(x,t)+\bA_\ep \bv^\ep(x,t)=\bdelta^\ep(x), & & x\in\bbR_\ep^d,\ t>0,\\
& \bv^\ep(x,0)=0,   & & x\in\bbR_\ep^d.
\end{aligned}
\right.
\end{equation}
Here $\dot{}=\partial/\partial t$.
\end{definition}

 In Section~\ref{secFirstGreen}, we will establish a
connection between the discrete and continuous Green functions in
terms of asymptotic formulas.

Using the discrete Fourier transform~\eqref{eqFourier}, its
inverse~\eqref{eqInverseFourier}, and
relation~\eqref{eqFourierSymbol}, we obtain the explicit
representations
\begin{align}
\bu^\ep(x,t)&=\dot \bv^\ep(x,t)
=\dfrac{1}{(2\pi)^d}\int_{R_{\pi\ep^{-1}}}
e^{-t\ep^{-\ell}A(\eta\ep)}e^{i
x\eta}\,d\eta,\label{eqz_nFourierMultiDim}\\
\bv^\ep(x,t)&=\dfrac{1}{(2\pi)^d}\int_{R_{\pi\ep^{-1}}}
\dfrac{1-e^{-t\ep^{-\ell}A(\eta\ep)}}{\ep^{-\ell}A(\eta\ep)}\,e^{i
x\eta}\,d\eta.\label{eqy_nFourierMultiDim}
\end{align}

Changing the variables in the integrals
in~\eqref{eqz_nFourierMultiDim} and~\eqref{eqy_nFourierMultiDim},
we obtain for $J=0,1,2,\dots$
\begin{equation}\label{eqRescalingGreen}
\dfrac{\partial ^J \bu^\ep(x,t)}{\partial t^J}\equiv
\ep^{-J\ell-d} \dfrac{\partial^J
\bu^1\left(x',\tau\right)}{\partial
\tau^J}\Bigg|_{x'=x/\ep,\,\tau=t/\ep^\ell},\quad
\bv^\ep(x,t)\equiv \ep^{\ell-d}
\bv^1\left(\dfrac{x}{\ep},\dfrac{t}{\ep^\ell}\right).
\end{equation}

\section{Asymptotics of the first Green function
$\bu^\ep(x,t)$}\label{secFirstGreen}

\subsection{Formulation of the result}

 In the theorems
below, we will not explicitly indicate the dependence of the
remainders in asymptotic formulas on   the number $M$ fixed in
Condition~\ref{condAppriximation}.

In this section, we obtain the following asymptotics for
$\bu^\ep(x,t)$.

\begin{theorem}\label{thAsymphMultiDim}
For any $\ep>0$, $t_0>0$, integer $K\ge M$, integer $J\ge 0$, and
all $x\in\bbR_\ep^d$ and $t\ge t_0\ep^\ell$, we have
\begin{equation}\label{eqAsymphGeneralMultiDim}
\dfrac{\partial^J\bu^\ep(x,t)}{\partial t^J}=\sum\limits_{k=0}^K
\dfrac{\ep^k}{t^{(k+d)/\ell+J}}\,
{H}_{Jk}\left(\dfrac{x}{t^{1/\ell}}\right) + \br_\bu^\ep(J,K; x,t),
\end{equation}
where
\begin{equation}\label{eqHk}
{H}_{Jk}(y):=\left\{
\begin{aligned}
&(-\cA(\cD))^J{H}(y) & &\text{if } k=0,\\
&0   & &\text{if } k=1,\dots,M-1,\\
&\text{finite linear combinations}\\
&\text{of derivatives of ${H}(y)$} & &\text{if } k=M,\dots,K,
\end{aligned}
\right.
\end{equation}
${H}(y)$ is given by~\eqref{eqhGeneral},
\begin{equation}\label{eqHkRemainder}
|\br_\bu^\ep(J,K; x,t)|\le
\dfrac{\ep^{K+1}R_\bu(J,K,t_0)}{t^{(K+d+1)/\ell+J}},
\end{equation}
and $R_\bu(J,K,t_0)\ge 0$ does not depend on $\ep>0$,
$x\in\bbR_\ep^d$, and $t\ge t_0\ep^\ell$.
\end{theorem}

\begin{remark}
It follows from the proof below that
${H}_{Jk}(y)=\cR_{J,k+\ell}(\cD){H}(y)$, where $\cR_{J,k+\ell}(\xi)$
are polynomials of degree $\ge k+\ell$ defined by~\eqref{eqAsymphWeak4'5MultiDimNew}.
\end{remark}

Note that the main term in the asymptotic representation of
$\bu^\ep(x,t)$ coincides with the continuous Green
function~\eqref{eqContGreen}.
 Theorem~\ref{thAsymphMultiDim}
implies that the higher the approximation order~$M$ (see
Definition~\ref{defAppriximation}) is, the faster the first
discrete Green function converges to the continuous Green
function, as $t\to\infty$ or as $\ep\to0$. More precisely, the
following corollaries hold.

\begin{corollary}\label{corApproxt} Under the assumptions of
Theorem~$\ref{thAsymphMultiDim}$, for any $\ep_0>0$,
$$
\sup\limits_{\ep\in(0,\ep_0]}\sup\limits_{x\in\bbR_\ep^d}
\left|\bu^\ep(x,t)-\dfrac{1}{t^{d/\ell}}\,
{H}\left(\dfrac{x}{t^{1/\ell}}\right)\right|=O\left(\dfrac{1}{t^{(M+d)/\ell}}\right)\quad\text{as
} t\to\infty.
$$
\end{corollary}

\begin{corollary}\label{corApproxEpsilon} Under the assumptions of
Theorem~$\ref{thAsymphMultiDim}$, for any $t_1>0$,
$$
\sup\limits_{t\ge t_1}\sup\limits_{x\in\bbR_\ep^d}
\left|\bu^\ep(x,t)-\dfrac{1}{t^{d/\ell}}\,
{H}\left(\dfrac{x}{t^{1/\ell}}\right)\right|=O(\ep^M)\quad\text{as
} \ep\to0.
$$
\end{corollary}

\subsection{Proof of Theorem~$\ref{thAsymphMultiDim}$}

{\bf Setp 1.} First, we assume that $\ep=1$ and denote
$\bu(x,t):=\bu^1(x,t)$. Note that $|\bu(x,t)|$ is bounded for
$x\in\bbZ^d$ and $t\ge 0$ due to~\eqref{eqz_nFourierMultiDim}
and~\eqref{eqDiscreteEllipticitym}. Therefore, if we show that
Theorem~$\ref{thAsymph}$ holds with $\ep=1$ and some $t_0>0$, it
will follow that the theorem holds with $\ep=1$ and all $t_0>0$.
Let us fix $t_0=1$ and consider $t\ge 1$.
%

Differentiating~\eqref{eqz_nFourierMultiDim} with $\ep=1$ and making the change
of variables $\xi=t^{1/\ell}\eta$, we
  write
\begin{equation}\label{eqAsymphWeak1Multidim}
\dfrac{\partial^J \bu(x,t)}{\partial t^J}  =\dfrac{1}{(2\pi)^d\,
t^{d/\ell}}\int\limits_{R_{\pi
 t^{1/\ell}}} (-1)^J \left(A(t^{-1/\ell}\xi)\right)^J
   e^{-tA(t^{-1/\ell}\xi)} \cdot e^{i
x t^{-1/\ell}\xi}\, d\xi=I_1(x,t)+I_2(x,t),
\end{equation}
where
\begin{equation}\label{eqAsymphWeak4Multidim}
I_1(x,t):=\dfrac{1}{(2\pi)^d\, t^{d/\ell}}\int\limits_{R_{\pi
 t^{1/L}}} (-1)^J \left(A(t^{-1/\ell}\xi)\right)^J
   e^{-tA(t^{-1/\ell}\xi)} \cdot e^{i
x t^{-1/\ell}\xi}\, d\xi,
\end{equation}
\begin{equation}\label{eqAsymphWeak4'Multidim}
I_2(x,t):=\dfrac{1}{(2\pi)^d\, t^{d/\ell}}\int\limits_{R_{\pi
 t^{1/\ell}}\setminus R_{\pi
 t^{1/L}}} (-1)^J \left(A(t^{-1/\ell}\xi)\right)^J
   e^{-tA(t^{-1/\ell}\xi)} \cdot e^{i
x t^{-1/\ell}\xi}\, d\xi,
\end{equation}
\begin{equation}\label{eqAsymphWeak4'MultidimL}
L:=(\ell+2)\ell.
\end{equation}

{\bf Step 2.} Let us estimate $I_1(x,t)$
in~\eqref{eqAsymphWeak4Multidim}. First, we consider the term
$e^{-tA(t^{-1/\ell}\xi)}$. By Lemma~\ref{lConnectionSymbols},
\begin{equation}\label{eqAsympAthetaMultiDim}
t A(t^{-1/\ell}\xi)=\cA(\xi)+\sum\limits_{k=M}^K
t^{-k/\ell}\cA_{k+\ell}(\xi)+\hat A(t,\xi),
\end{equation}
where
\begin{equation}\label{eqAsympAthetaMultiDim'}
|\hat A(t,\xi)|\le c_1 t^{-(K+1)/\ell} |\xi|^{K+1+\ell},\quad
\xi\in R_{\pi t^{1/\ell}},\ t\ge 1.
\end{equation}
Here and below $c_1,c_2,{\dots}>0$ do not depending on $\xi$ and
$t$.

Further,  \eqref{eqAsymphWeak4'MultidimL} implies for $\xi\in
R_{\pi t^{1/L}}$ and $t\ge 1$ that
$$
t^{-k/\ell}|\xi|^{k+\ell}\le c_2 t^{-k/\ell}t^{(k+\ell)/L}= c_2
t^{(\ell-(\ell+1) k)/L} \le c_2,\quad k=1,2,\dots.
$$
Therefore,
\begin{equation}\label{eqAsymphWeak4'3MultiDim}
\left|\sum\limits_{k=M}^K t^{-k/\ell}\cA_{k+\ell}(\xi)+\hat
A(t,\xi)\right|\le c_3, \quad \xi\in R_{\pi t^{1/L}},\ t\ge 1.
\end{equation}

Using the Taylor expansion of the exponent function, we obtain
\begin{equation}\label{eqAsymphWeak4'5MultiDim}
\begin{aligned}
 \exp\left(-\sum\limits_{k=M}^K
t^{-k/\ell}\cA_{k+\ell}(\xi) - \hat A(t,\xi)\right) &
=1+\sum\limits_{m=1}^K \dfrac{1}{m!}\left(\sum\limits_{k=M}^K
t^{-k/\ell}\cA_{k+\ell}(\xi)+\hat A(\xi)\right)^m + \hat A_1(\xi,t)\\
  &= 1+ \sum\limits_{k=M}^K t^{-k/\ell}\cP_{k+\ell}(\xi)  +
 \hat P(\xi,t)
\end{aligned}
\end{equation}
Here $\cP_{k+\ell}(\xi)$ are polynomials of degree $\ge k+\ell$ and the remainders are estimated due to
relations~\eqref{eqAsympAthetaMultiDim'}
and~\eqref{eqAsymphWeak4'3MultiDim} as follows:
$$
|\hat A_1(\xi,t)|\le e^{c_3}\left|\sum\limits_{k=M}^K
t^{-k/\ell}\cA_{k+\ell}(\xi)+\hat A(\xi,t)\right|^{K+1}\le c_4
e^{c |\xi|^\ell/2} t^{-M(K+1)/\ell}, \quad \xi\in R_{\pi
t^{1/L}},\ t\ge 1,
$$
\begin{equation}\label{eqAsymphWeak4'5MultiDim'}
|\hat P(\xi,t)|\le  c_5 e^{c |\xi|^\ell/2} t^{-(K+1)/\ell}, \quad
\xi\in R_{\pi t^{1/L}},\ t\ge 1,
\end{equation}
where  $c$ is the constant from~\eqref{eqcAelliptic}.

Now we estimate the term $\left(A(t^{-1/\ell}\xi)\right)^J$.
Using~\eqref{eqAsympAthetaMultiDim}, we have
\begin{equation}\label{eqAsymphWeak4'5MultiDimJ}
\left(A(t^{-1/\ell}\xi)\right)^J=t^{-J}\left((\cA(\xi))^J+\sum\limits_{k=M}^K
t^{-k/\ell} \cD_{J,k+\ell}(\xi)+\hat D(t,\xi)\right),
\end{equation}
where $\cD_{J,k+\ell}(\xi)$ are polynomials of degree $\ge k+\ell$
and
\begin{equation}\label{eqAsymphWeak4'5MultiDim'J}
|\hat D(\xi,t)|\le  c_6 e^{c |\xi|^\ell/2} t^{-(K+1)/\ell}, \quad
\xi\in R_{\pi t^{1/\ell}},\ t\ge 1.
\end{equation}

Combining~\eqref{eqAsymphWeak4'5MultiDim}--\eqref{eqAsymphWeak4'5MultiDim'J}
with~\eqref{eqAsympAthetaMultiDim} yields
\begin{equation}\label{eqAsymphWeak4'5MultiDimNew}
\begin{aligned}
&(-1)^J \left(A(t^{-1/\ell}\xi)\right)^J
   e^{-tA(t^{-1/\ell}\xi)}\\
   &\quad =e^{-\cA(\xi)}\left(   t^{-J}
   \left(-\cA(\xi)\right)^J + \sum\limits_{k=M}^K
t^{-k/\ell-J} \cR_{J,k+\ell}(\xi)+\hat R(t,\xi)\right),
\end{aligned}
\end{equation}
where $\cR_{J,k+\ell}(\xi)$ are polynomials of degree $\ge k+\ell$
and
\begin{equation}\label{eqAsymphWeak4'5MultiDimNewRemainder}
|\hat R(\xi,t)|\le  c_7 e^{c |\xi|^\ell/2} t^{-(K+1)/\ell-J},
\quad \xi\in R_{\pi t^{1/L}},\ t\ge 1.
\end{equation}

Substituting~\eqref{eqAsymphWeak4'5MultiDimNew}
into~\eqref{eqAsymphWeak4Multidim}, we have
\begin{equation}\label{eqAsymphWeak4'7MultiDim}
\begin{aligned}
 I_1(x,t) & = \dfrac{1}{(2\pi)^d\, t^{d/\ell+J}} \int\limits_{R_{\pi
t^{1/L}}}
   \left(-\cA(\xi)\right)^J e^{-\cA(\xi)} \cdot e^{i
x t^{-1/\ell}\xi}\, d\xi\\
& + \sum\limits_{k=M}^K \dfrac{1}{(2\pi)^d\, t^{(k+d)/\ell+J}}
\int\limits_{R_{\pi t^{1/L}}} \cR_{J,k+\ell}(\xi)\,e^{-\cA(\xi)}
\cdot e^{i x t^{-1/\ell}\xi}\, d\xi \\
& + \dfrac{1}{(2\pi)^d\, t^{d/\ell}} \int\limits_{R_{\pi t^{1/L}}}
e^{-\cA(\xi)}\hat R(\xi,t) \cdot e^{i x t^{-1/\ell}\xi}\,d\xi.
\end{aligned}
\end{equation}

Note that, due to~\eqref{eqcAelliptic}, the polynomials
$\cB(\xi)=(-\cA(\xi))^J$ and $\cB(\xi)=\cR_{J,k+l}(\xi)$ satisfy
\begin{equation}\label{eqAsymphWeak4'7MultiDim'}
\begin{aligned}
\left|\,\int\limits_{\bbR^d\setminus R_{\pi t^{1/L}}}
\cB(\xi) e^{-\cA(\xi)}  \cdot e^{i x t^{-1/\ell}\xi}\, d\xi\right|&
\le
  \int\limits_{\bbR^d\setminus R_{\pi t^{1/L}}}  |\cB(\xi)|\cdot e^{-c
|\xi|^\ell} \, d\xi =
O\left(\dfrac{1}{t^{(K+1)/\ell}}\right),
\end{aligned}
\end{equation}
 where $O(\cdot)$ is taken as $t\to\infty$ (uniformly with
respect to $x$). Similarly, due
to~\eqref{eqAsymphWeak4'5MultiDimNewRemainder}, $\hat R(\xi,t)$ satisfies
\begin{equation}\label{eqAsymphWeak4'7MultiDim''}
\begin{aligned}
\left|\,\int\limits_{ R_{\pi t^{1/L}}} e^{-\cA(\xi)}
\hat R(\xi,t) \cdot e^{i x t^{-1/\ell}\xi}\, d\xi\right| & \le
\dfrac{c_7}{t^{(K+1)/\ell+J}} \int\limits_{ R_{\pi t^{1/L}}}
e^{-\Re
\cA(\xi)}\cdot e^{c |\xi|^\ell/2} \, d\xi \\
& \le \dfrac{c_7}{t^{(K+1)/\ell+J}} \int\limits_{\bbR^d} e^{-c
|\xi|^\ell/2} \, d\xi.
\end{aligned}
\end{equation}
Combining~\eqref{eqAsymphWeak4'7MultiDim}--\eqref{eqAsymphWeak4'7MultiDim''},
we obtain
\begin{equation}\label{eqAsymphWeak4'5MultiDim'''}
\begin{aligned}
 I_1(x,t) & = \dfrac{1}{(2\pi)^d\, t^{d/\ell+J}}
 \int\limits_{\bbR^d}
(-\cA(\xi))^J e^{-\cA(\xi)} \cdot e^{i
x t^{-1/\ell}\xi}\, d\xi\\
& + \sum\limits_{k=M}^K \dfrac{1}{(2\pi)^d\, t^{(k+d)/\ell+J}}
\int\limits_{\bbR^d} \cR_{J,k+\ell}(\xi)\,e^{-\cA(\xi)} \cdot e^{i x
t^{-1/\ell}\xi}\, d\xi +
O\left(\dfrac{1}{t^{(K+d+1)/\ell+J}}\right).
\end{aligned}
\end{equation}

 Relations~\eqref{eqhGeneral} and
\eqref{eqAsymphWeak4'5MultiDim'''} and the fact that
$\cR_{J,k+\ell}(\xi)$ are polynomials imply
\begin{equation}\label{eqAsymphWeak4'8MultiDim}
I_1(x,t)=\dfrac{1}{t^{d/\ell+J}}\,{H}_{J0}\left(\dfrac{x}{t^{1/\ell}}\right)+
\sum\limits_{k=M}^K\dfrac{1}{t^{(k+d)/\ell+J}}\,
{H}_{Jk}\left(\dfrac{x}{t^{1/\ell}}\right)+O\left(\dfrac{1}{t^{(K+d+1)/\ell+J}}\right)
\end{equation}
where ${H}_{J0}(y)=(-\cA(\cD))^J{H}(y)$ and
${H}_{Jk}(y)=\cR_{J,k+\ell}(\cD){H}(y)$.

{\bf Step 3.} Let us estimate $I_2(x,t)$
in~\eqref{eqAsymphWeak4'Multidim}.
Using~\eqref{eqDiscreteEllipticitymUpperBound} and
\eqref{eqDiscreteEllipticitym}, we have for $\xi\in R_{\pi
 t^{1/\ell}}\setminus R_{\pi
 t^{1/L}}$
$$
\left|\left(A(t^{-1/\ell}\xi)\right)^J e^{-t
A(t^{-1/\ell}\xi)}\right|\le c_8 e^{-c|\xi|^\ell/2}\le c_8 e^{-c
\pi^\ell \, t^{\ell/L}/2}.
$$
Therefore,
\begin{equation}\label{eqAsymphWeak7MultiDim}
|I_2(x,t)|\le c_9 e^{-c \pi^\ell \, t^{\ell/L}/2},\quad t\ge 1.
\end{equation}

Relations~\eqref{eqAsymphWeak1Multidim}--\eqref{eqAsymphWeak4'Multidim},
\eqref{eqAsymphWeak4'8MultiDim} and~\eqref{eqAsymphWeak7MultiDim}
yield the assertion of Theorem~\ref{thAsymphMultiDim} with
$\ep=1$. Using the first relation in~\eqref{eqRescalingGreen}, we
obtain the assertion with any $\ep>0$.

\section{Preliminary asymptotics of the second Green function
$\bv^\ep(x,t)$}\label{secSecondGreenPrelim}

Throughout this section, we assume that $\ep=1$ and denote
$\bv(x,t):=\bv^1(x,t)$. First, we obtain an asymptotic formula for
$\bv(x,t)$ up to some unknown  function $\Omega(x)$. In the
next   section, we will provide an explicit formula for it.

We introduce the spherical coordinates $(r,\ph)$ with $r>0$ and
$\ph\in\Phi$, where
\begin{equation}\label{eqSphericalCoordinatesPhi}
\Phi:=\{\ph\in\bbR^{d-1}: \ph_1,\dots,\ph_{d-2}\in[0,\pi],\
\ph_{d-1}\in[0,2\pi)\}.
\end{equation}
They are related to the Cartesian coordinates $y\in\bbR^d$ via
\begin{equation}\label{eqSphericalCoordinates}
\begin{aligned}
&y_1=r\cos\ph_1,\quad y_2=r\sin \ph_1 \cos\ph_2,\quad  y_3=r\sin
\ph_1
\sin\ph_2 \cos\ph_3,\dots,\\
& y_{d-1}=r\sin \ph_1\cdots \sin\ph_{d-2} \cos\ph_{d-1},\quad
y_d=r\sin \ph_1\cdots \sin\ph_{d-2} \sin\ph_{d-1},
\end{aligned}
\end{equation}
and the Jacobian is given by
\begin{equation}\label{eqSphericalCoordinatesJacobian}
dy=r^{d-1}{J}(\ph)\,dr\, d\ph,\quad {J}(\ph):=\sin^{d-2}\ph_1
\sin^{d-3} \ph_2\cdots\sin\ph_{d-2}.
\end{equation}

We will denote functions written in Cartesian and spherical
coordinates by the same letter. For example ${H}(r,\ph)$ stands
for ${H}(y)$, etc.

We set
$$
{H}_k(y):={H}_{0k}(y),\quad y\in\bbR^d,
$$
where ${H}_{0k}(y)$ are the functions from
Theorem~\ref{thAsymphMultiDim} with $J=0$.
 For $y\in\bbR^d\setminus\{0\}$, we introduce the functions
\begin{align}
{F}_k(y):=
 \ell r^{\ell-d-k} & \int_r^\infty
\dfrac{{H}_k(\rho,\ph)}{\rho^{\ell-d-k+1}}\,d\rho & &\text{if}\
  k=0,\dots,\ell-d, \label{eqcFLog-Corrected1}\\
{F}_k(y):= -\dfrac{\ell}{r^{d+k-\ell}} & \int_0^r
  \rho^{d+k-\ell-1}{H}_k(\rho,\ph)\,d\rho & &\text{if}\
k\ge \ell-d+1. \label{eqcFLog-Corrected2}
\end{align}
Definition~\eqref{eqcFLog-Corrected1} will be used only in the
case $\ell-d\ge 0$ and definition~\eqref{eqcFLog-Corrected2} only
for $k\ge0$.

\begin{remark}\label{remcFkEqual0}
Due to~\eqref{eqHk}, ${F}_k(y)\equiv 0$ for $k=1,\dots,M-1$.
\end{remark}

\begin{lemma}\label{lFk12} Let  a function ${F}_k(y)$ be defined
by~\eqref{eqcFLog-Corrected1} or~\eqref{eqcFLog-Corrected2}. Then
\begin{enumerate}

\item\label{lFk2} it  satisfies the partial differential equation
\begin{equation}\label{eqEquationfkMultiDim}
(\ell-d-k){F}_k-\sum\limits_{j=1}^d
y_j\dfrac{\partial{F}_k}{\partial y_j} =\ell{H}_k,\quad
y\in\bbR^d\setminus\{0\}.
\end{equation}

\item\label{lFk1} it tends to $0$ as $|y|\to\infty$,

\item\label{lFk3} if $k\ne \ell-d$, it is bounded in $\bbR^d$ and,
for each fixed $\ph$, continuous with respect to $r$ at the
origin,

\item\label{lFk4} if $k\ge \ell-d+1$, it is continuous with
respect to $y$ at the origin and
\begin{equation}\label{eqFk4}
{F}_k(0):=\lim\limits_{y\to
0}{F}_k(y)=-\dfrac{\ell}{k+d-\ell}{H}_k(0).
\end{equation}

\end{enumerate}
\end{lemma}
\proof  Item~\ref{lFk2} follows by observing that the equation
in~\eqref{eqEquationfkMultiDim} is equivalent to the following:
$$
(\ell-d-k){F}_k-r\dfrac{\partial{F}_k}{\partial r}=\ell{H}_k,\quad
r>0.
$$
The latter can be solved as an ordinary differential equation for
each fixed $\ph\in\Phi$. It is easy to check that the unique
solution that is bounded at infinity is ${F}_k(y)$ defined
by~\eqref{eqcFLog-Corrected1} or~\eqref{eqcFLog-Corrected2},
respectively.

Items~\ref{lFk1}, \ref{lFk3}, and~\ref{lFk4} follow from the fact
that ${H}_k(y)$ are infinitely differentiable at the origin and
are rapidly decreasing functions.
\endproof

Now we find a preliminary asymptotic formula for $\bv(x,t) $ in
terms of the above functions ${F}_k(y)$. It contains an unknown
grid function $\Omega(x)$ and, if $\ell-d\ge 0$, it is not valid
at $x=0$. The reason for the latter is essentially the fact that
the functions ${F}_k(y)$ are not continuous at $y=0$ for
$k=0,\dots, \ell-d$.

\begin{lemma}\label{lAsymphMultiDimPrelim}
 There exists a grid function $\Omega(x)$ such that,
for any $t_0>0$, integer $K\ge \max(M,\ell-d)$,  $x\in
\bbZ^d\setminus\{0\}$, and  $t\ge t_0$,
\begin{equation}\label{eqAsympfGeneralPrelimLess0}
\bv(x,t)  = \Omega(x) + \sum\limits_{k=0}^K
\dfrac{1}{t^{(k+d)/\ell-1}}\,
{F}_k\left(\dfrac{x}{t^{1/\ell}}\right)  + \br_\bv(K; x,t).
\end{equation}
 Here ${F}_k(y)$ are given by~\eqref{eqcFLog-Corrected1}
and~\eqref{eqcFLog-Corrected2},
\begin{equation}\label{eqFkRemainderPrelim}
|\br_\bv(K; x,t)|\le \dfrac{R_\bv(K,t_0)}{t^{(K+d+1)/\ell-1}},
\end{equation}
and $R_\bv(K,t_0)\ge 0$ does not depend on $x\in\bbZ^d$ and $t\ge
t_0$.

Furthermore, if $\ell-d\le -1$, then the equality
in~\eqref{eqAsympfGeneralPrelimLess0} and the estimate
in~\eqref{eqFkRemainderPrelim} also hold  for $x=0$ with
${F}_k(0)$ defined by~\eqref{eqFk4}.
\end{lemma}
\proof We integrate~\eqref{eqAsymphGeneralMultiDim} with $J=0$ and
$\ep=1$ from $t_0$ to $t$ and use the identity
$\bu(x,t)\equiv\dot\bv(x,t)$:
\begin{equation}\label{eqAsympf1MultiDim}
\begin{aligned}
 \bv(x,t)-\bv(x,t_0) & =  \sum\limits_{k=0}^K \int\limits_{t_0}^t
\dfrac{1}{\tau^{(k+d)/\ell}}\,
{H}_k\left(\dfrac{x}{\tau^{1/\ell}}\right)\,d\tau\\
& +\int\limits_{t_0}^\infty \br_\bu^1(0,K; x,\tau)\,
d\tau-\int\limits_{t}^\infty \br_\bu^1(0,K; x,\tau)\,
d\tau,\quad x\in\bbZ^d.
\end{aligned}
\end{equation}

It follows from the chain rule and~\eqref{eqEquationfkMultiDim}
that
\begin{equation}\label{eqAsympf2MultiDim}
\dfrac{1}{\tau^{(k+d)/\ell}}\,
{H}_k\left(\dfrac{x}{\tau^{1/\ell}}\right)=\dfrac{\partial}{\partial\tau}\left(\dfrac{1}{\tau^{(k+d)/\ell-1}}\,
{F}_k\left(\dfrac{x}{\tau^{1/\ell}}\right)\right),\quad
x\in\bbZ^d\setminus\{0\}.
\end{equation}

Combining~\eqref{eqAsympf1MultiDim} and~\eqref{eqAsympf2MultiDim}
and taking into account~\eqref{eqHkRemainder} and the inequality
$K\ge \ell-d$ yields
\begin{equation}\label{eqAsympf4MultiDim}
\bv(x,t)  = \Omega(x) + \sum\limits_{k=0}^K
\dfrac{1}{t^{(k+d)/\ell-1}}\,
{F}_k\left(\dfrac{x}{t^{1/\ell}}\right)  + \br_\bv(K; x,t),\quad x\in\bbZ^d\setminus\{0\},
\end{equation}
with some grid function $\Omega(x)$ and the remainder
$\br_\bv(K; x,t)$ satisfying~\eqref{eqFkRemainderPrelim}. This
proves~\eqref{eqAsympfGeneralPrelimLess0} in the case
 $x\ne0$. If $\ell-d\le -1$, the same formula for $x=0$ follows
directly from~\eqref{eqAsympf1MultiDim} by substituting $x=0$
therein and using~\eqref{eqFk4}.

To see that $\Omega(x)$ and $\br_\bv(K; x,t)$ do not depend on $t_0$, denote by
$\Omega^*(x)$ and $\br_\bv^*(K; x,t)$ the functions corresponding to $t_0$ replaced by some
$t_0^*>t_0$. Then~\eqref{eqAsympf4MultiDim} implies that
$$
0=\Omega(x)-\Omega^*(x)+\br_\bv(K; x,t)-\br_\bv^*(K;
x,t),\quad t\ge t_0^*.
$$
Passing to the limit as $t\to\infty$ yields
$\Omega(x)=\Omega^*(x)$. In turn, this and~\eqref{eqAsympf4MultiDim} imply that  $\br_\bv(K; x,t)=\br_\bv^*(K; x,t)$ for $t\ge t_0^*$.
\endproof

In the case~$\ell-d\ge 0$, we will need modified  versions of the
functions in~\eqref{eqcFLog-Corrected1}. To introduce them, we
denote by $n_\ph$ the unit vector originated at $0$ and pointing
in the direction $\ph$ and by $f^{(j)}(0,\ph)$ the $j$th
directional derivative of a function $f(y)$ in the direction
$n_\ph$ at the origin. Now, for each fixed $\ph\in\Phi$, we expand
${H}_k(\rho,\ph)$ into the Taylor series around $\rho=0$:
\begin{equation}\label{eqFindGamma2'}
{H}_k(\rho,\ph)=\sum\limits_{j=0}^{\ell-d-k}
\dfrac{{H}_k^{(j)}(0,\ph)}{j!}\rho^j + \hat{H}_k(\rho,\ph),\quad
k=0,\dots,\ell-d,
\end{equation}
where
\begin{equation}\label{eqFindGamma3}
|\hat{H}_k(\rho,\ph)|\le c_1 \rho^{\ell-d-k+1},\quad \rho\in
[0,1],
\end{equation}
and $c_1>0$ does not depend on $\rho\in[0,1]$ and $\ph\in\Phi$.

For $\ell-d\ge 0$ and $y>0$, we set
%
\begin{equation}\label{eqFormulafkMultiDim1}
\hat{F}_k(y):=\ell r^{\ell-d-k}\int\limits_{r}^1
\dfrac{\hat{H}_k(\rho,\ph)}{\rho^{\ell-d-k+1}}\,d\rho,\quad
k=0,\dots,\ell-d,
\end{equation}
\begin{equation}\label{eqFindGamma5'}
 \hat{F}(y):=\ell\sum\limits_{k=0}^{\ell-d}
r^{\ell-d-k}\int\limits_{0}^1
\dfrac{\hat{H}_k(\rho,\ph)}{\rho^{\ell-d-k+1}}\,d\rho,
\end{equation}
\begin{equation}\label{eqcGk}
 {G}_k(y):=\left\{
\begin{aligned}
&\dfrac{\ell}{\ell-k-d} \sum\limits_{m=0}^{k}
\dfrac{{H}_m^{(k-m)}(0,\ph)\, r^{k-m}}{(k-m)!} & & \text{if }
k=0,\dots,\ell-d-1,\\
& \sum\limits_{m=0}^{\ell-d} \dfrac{{H}_m^{(\ell-d-m)}(0,\ph) \,
r^{\ell-d-m}}{(\ell-d-m)!} & & \text{if } k=\ell-d,
\end{aligned}
\right.
\end{equation}
\begin{equation}\label{eqFindGamma4'}
\begin{aligned}
 \Omega_k(y)&:=\ell r^{\ell-d-k} \Bigg(  \int\limits_1^\infty
\dfrac{{H}_k(\rho,\ph)}{\rho^{\ell-d-k+1}}\,d\rho \\
 & -   \sum\limits_{j=0}^{\ell-d-k-1}
\dfrac{{H}_k^{(j)}(0,\ph)}{j!(\ell-j-d-k)} -
\dfrac{{H}_k^{(\ell-d-k)}(0,\ph)}{(\ell-d-k)!}\ln r\Bigg),\quad
k=0,\dots,\ell-d.
\end{aligned}
\end{equation}

\begin{remark}
\begin{enumerate}
\item The integral in~\eqref{eqFindGamma5'} converges due
to~\eqref{eqFindGamma3}.

\item The first formula in~\eqref{eqcGk} is used only if
$\ell-d\ge 1$.

\item In~\eqref{eqFindGamma4'} and below, we use the convention
that the sum vanishes if the upper summation limit is less than
the lower limit.
\end{enumerate}
\end{remark}

The following assertion follows from~\eqref{eqFindGamma3} and the
definition of $\hat{F}_k(y)$, $\hat{F}(y)$, and ${G}_k(y)$.

\begin{lemma}\label{lFk12Modified} Let $\ell-d\ge 0$.
\begin{enumerate}
\item The functions $\hat{F}_k(y)$, $\hat{F}(y)$, and ${G}_k(y)$
are bounded near the origin. For each fixed $\ph\in\Phi$, they are
continuous with respect to $r\ge0$.

\item\label{lFk12Modified2} We have
$$
\lim\limits_{t\to\infty}\sum\limits_{k=0}^{\ell-d}
t^{1-(d+k)/\ell}\hat{F}_k\left(\dfrac{x}{t^{1/\ell}}\right)=\hat{F}(x),\quad
x\in\bbZ^d\setminus\{0\}.
$$
\end{enumerate}
\end{lemma}

Now, based on Lemma~\ref{lAsymphMultiDimPrelim} in the case
$\ell-d\ge 0$, we find another preliminary asymptotic formula for
$\bv(x,t) $ in terms of the above functions ${F}_k(y)$,
$\hat{F}_k(y)$, $\hat{F}(y)$, ${G}_k(y)$, and $\Omega_k(y)$. Now
the asymptotic formula is extended to $x=0$, but still contains an
unknown grid function $\Omega(x)$, which will be determined in the
next  section.

\begin{lemma}\label{lAsymphMultiDimPrelimVer2}
Let $\ell-d\ge 0$. Then, for any $t_0>0$, integer $K\ge
\max(M,\ell-d)$,   and  $t\ge t_0$, the following hold.
\begin{enumerate}
\item If $x\in \bbZ^d\setminus\{0\}$, then, along with
asymptotics~\eqref{eqAsympfGeneralPrelimLess0}, we have
\begin{equation}\label{eqFindGamma5}
\begin{aligned}
\bv(x,t) & = \sum\limits_{k=0}^{\ell-d-1}
t^{1-(d+k)/\ell}\left(\hat{F}_k\left(\dfrac{x}{t^{1/\ell}}\right)
+ {G}_k(x)\right)
+  {G}_{\ell-d}(x)\ln t\\
& +  \hat{F}_{\ell-d}\left(\dfrac{x}{t^{1/\ell}}\right) +
\hat\Omega(x) + \sum\limits_{k=\ell-d+1}^K
\dfrac{1}{t^{(d+k)/\ell-1}}
 {F}_k\left(\dfrac{x}{t^{1/\ell}}\right)+ \br_\bv(K; x,t),
\end{aligned}
\end{equation}
where
\begin{equation}\label{eqhatGamma}
\hat\Omega(y)=\Omega(y) + \sum\limits_{k=0}^{\ell-d}\Omega_k(y).
\end{equation}

\item\label{eqFindGamma5Part2} If $x=0$, then we have
\begin{equation}\label{eqFindGamma5X0}
\begin{aligned}
\bv(0,t) & = \sum\limits_{k=0}^{\ell-d-1}  t^{1-(d+k)/\ell}\dfrac{\ell\,
{H}_k(0)}{\ell-d-k} + {H}_{\ell-d}(0)\ln t\\
&  + \omega + \sum\limits_{k=\ell-d+1}^K
\dfrac{1}{t^{(d+k)/\ell-1}}  \dfrac{\ell{H}_k(0)}{\ell-d-k}+ \br_\bv(K; 0,t),
\end{aligned}
\end{equation}
where $\omega\in\bbC$  does not depend on $t_0$ and $t$ and $ {F}_k(0)$ are defined
in~\eqref{eqFk4}.
\end{enumerate}
In both cases, $\br_\bv(K; x,t)$ satisfies
estimate~\eqref{eqFkRemainderPrelim}.
\end{lemma}
\proof First, we fix $x\in\bbZ^d\setminus\{0\}$. Consider the
terms in~\eqref{eqAsympfGeneralPrelimLess0} corresponding to
$k=0,\dots,\ell-d$. As we will see below, they (in general) tend
to infinity as $t\to\infty$. For each $m\in\{0,\dots,\ell-d\}$, we
have
\begin{equation}\label{eqFindGamma2}
\begin{aligned}
t^{1-(m+d)/\ell} {F}_m\left(\dfrac{x}{t^{1/\ell}}\right) & =  \ell
r^{\ell-d-m}   \int\limits_{r t^{-1/\ell}}^1
\dfrac{{H}_m(\rho,\ph)}{\rho^{\ell-d-m+1}}\,d\rho + \ell
r^{\ell-d-m}   \int\limits_1^\infty
\dfrac{{H}_m(\rho,\ph)}{\rho^{\ell-d-m+1}}\,d\rho.
\end{aligned}
\end{equation}

Consider the first integral in~\eqref{eqFindGamma2}.
Using~\eqref{eqFindGamma2'}, we represent it as follows:
\begin{equation*}
\begin{aligned}
   \int\limits_{r t^{-1/\ell}}^1
\dfrac{{H}_m(\rho,\ph)}{\rho^{\ell-d-m+1}}\,d\rho=\sum\limits_{j=0}^{\ell-d-m}
\dfrac{{H}_m^{(j)}(0,\ph)}{j!} \int\limits_{r t^{-1/\ell}}^1
\rho^{j+d+m-\ell-1} \,d\rho + \int\limits_{r t^{-1/\ell}}^1
\dfrac{\hat{H}_m(\rho,\ph)}{\rho^{\ell-d-m+1}}\,d\rho,
\end{aligned}
\end{equation*}
which yields
$$
\begin{aligned}
&    \int\limits_{r t^{-1/\ell}}^1
\dfrac{{H}_m(\rho,\ph)}{\rho^{\ell-d-m+1}}\,d\rho  =
\sum\limits_{j=0}^{\ell-d-m-1} \dfrac{{H}_m^{(j)}(0,\ph)
r^{j+d+m-\ell}}{j!(\ell-j-d-m)}\, t^{1-(j+d+m)/\ell}
+\dfrac{{H}_m^{(\ell-d-m)}(0,\ph)}{ (\ell-d-m)!}\cdot\dfrac{\ln t}{\ell}\\
& \qquad - \sum\limits_{j=0}^{\ell-d-m-1}
\dfrac{{H}_m^{(j)}(0,\ph)}{j!(\ell-j-d-m)}
-\dfrac{{H}_m^{(\ell-d-m)}(0,\ph)}{(\ell-d-m)!}\ln r
  +
\int\limits_{r t^{-1/\ell}}^1
\dfrac{\hat{H}_m(\rho,\ph)}{\rho^{\ell-d-m+1}}\,d\rho.
\end{aligned}
$$
Combining this with~\eqref{eqFindGamma2}, we obtain
\begin{equation}\label{eqFindGamma4}
\begin{aligned}
 t^{1-(m+d)/\ell} {F}_m\left(\dfrac{x}{t^{1/\ell}}\right) & =
\ell \sum\limits_{j=0}^{\ell-d-m-1} \dfrac{{H}_m^{(j)}(0,\ph)\,
r^j}{j!(\ell-j-d-m)}\,
t^{1-(d+m+j)/\ell}\\
& +
\dfrac{{H}_m^{(\ell-d-m)}(0,\ph)\,r^{\ell-d-m}}{(\ell-d-m)!}\ln t
+\Omega_m(x)
  + t^{1-(m+d)/\ell}\hat{F}_m\left(\dfrac{x}{t^{1/\ell}}\right),
\end{aligned}
\end{equation}
where $\hat{F}_m(y)$ are given by~\eqref{eqFormulafkMultiDim1} and
$\Omega_m(y)$ by~\eqref{eqFindGamma4'}.

Equalities~\eqref{eqAsympfGeneralPrelimLess0}
and~\eqref{eqFindGamma4} imply~\eqref{eqFindGamma5} with
${G}_k(y)$ given by~\eqref{eqcGk} and $\hat\Omega(y)$
by~\eqref{eqhatGamma}.

For $x=0$, formula~\eqref{eqFindGamma5X0} directly follows  from
integrating~\eqref{eqAsympf1MultiDim} with $x=0$.
\endproof

The advantage of the representation~\eqref{eqFindGamma5} is that
it allows one to distinguish time-independent coefficients at
positive powers of $t$ in the asymptotic formula for $\bv(x,t)$.
Namely, the following assertion holds.

\begin{corollary}\label{corAsymphMultiDimPrelimVer2}
Let the assumptions of Lemma~$\ref{lAsymphMultiDimPrelimVer2}$
hold. Then, for each $x\in\bbZ^d\setminus\{0\}$,
\begin{equation}\label{eqAsymphMultiDimPrelimVer2PurePowers}
\begin{aligned}
\bv(x,t) & = \sum\limits_{k=0}^{\ell-d-1} t^{1-(d+k)/\ell}
{G}_k(x) +  {G}_{l-d}(x)\ln t +   \hat{F}(x) + \hat\Omega(x) +
o(1)\quad\text{as } t\to\infty,
\end{aligned}
\end{equation}
where $G_k(x)$, $\hat F(x)$, and $\hat\Omega(x)$ are given by~\eqref{eqcGk}, \eqref{eqFindGamma5'}, and~\eqref{eqhatGamma}, respectively,  and $o(1)$ depends on $x$.
\end{corollary}
\proof The assertion follows from~\eqref{eqFindGamma5} and Lemma~\ref{lFk12Modified} (part~\ref{lFk12Modified2}).
\endproof

\begin{remark} Asymptotic formulas~\eqref{eqAsympfGeneralPrelimLess0}
and~\eqref{eqFindGamma5} are uniform with respect to
$x\in\bbZ^d\setminus\{0\}$, in the sense that the constant
$R_\bv(K,t_0)$ in the estimate~\eqref{eqFkRemainderPrelim} does
not depend on $x\in\bbZ^d\setminus\{0\}$. Unlike those formulas,
the function $o(1)$ tends to $0$ as $t\to\infty$ in
relation~\eqref{eqAsymphMultiDimPrelimVer2PurePowers} uniformly
with respect to $x\ne0$ from any compact set, but {\em not}
uniformly in $\bbZ^d\setminus\{0\}$.

However, in the next subsection, we will find an asymptotics
analogous to~\eqref{eqAsymphMultiDimPrelimVer2PurePowers} directly
from the integral representation~\eqref{eqy_nFourierMultiDim} of
the Green function. In this new asymptotics, the time-independent
term will be explicitly present and thus will coincide with
$\hat{F}(x) + \hat\Omega(x)$
from~\eqref{eqAsymphMultiDimPrelimVer2PurePowers}. Hence, we will
find $\hat\Omega(x)$ and thus $\Omega(x)$ and ``finalize'' the
(uniform) asymptotic formulas~\eqref{eqAsympfGeneralPrelimLess0}
and~\eqref{eqFindGamma5}.
\end{remark}

\section{Asymptotics of the second Green function
$\bv^\ep(x,t)$}\label{secSecondGreen}

\subsection{Case $\ell-d\le -1$}\label{subsecEll-dNegative}

First, we complete Lemma~\ref{lAsymphMultiDimPrelim} by finding
$\Omega(x)$ in the case $\ell-d\le -1$.

\begin{theorem}\label{thAsympfGeneral1}
Let $\ell-d\le -1$. Then, for any $\ep>0$, $t_0>0$, integer $ K\ge
M$, and all $x\in\bbR_\ep^d$ and $t\ge t_0\ep^\ell$,
\begin{equation}\label{eqAsympfGeneralMultiDim-d-lBigger1}
\begin{aligned}
\bv^\ep(x,t)&=
\dfrac{1}{\ep^{d-\ell}}\Omega\left(\dfrac{x}{\ep}\right)+
\dfrac{1}{t^{d/\ell-1}} \,
{F}_0\left(\dfrac{x}{t^{1/\ell}}\right)\\
&+\sum\limits_{k=M}^K\dfrac{\ep^k}{t^{(k+d)/\ell-1}}\,
{F}_k\left(\dfrac{x}{t^{1/\ell}}\right) + \br_\bv^\ep(K; x,t),
\end{aligned}
\end{equation}
where
\begin{equation}\label{eqGammaldLess0}
\Omega(y):=\dfrac{1}{(2\pi)^d}\int_{R_\pi}\dfrac{e^{i
y\theta}}{A(\theta)}\,d\theta,
\end{equation}
 ${F}_k(y)$ are given by~\eqref{eqcFLog-Corrected2},
$$
|\br_\bv^\ep(K; x,t)|\le
\dfrac{\ep^{K+1}R_\bv(K,t_0)}{t^{(K+d+1)/\ell-1}},
$$
and $R_\bv(K,t_0)\ge 0$ does not depend on $\ep>0$,
$x\in\bbR_\ep^d$, and $t\ge t_0\ep^\ell$.
\end{theorem}
\proof Due to~\eqref{eqRescalingGreen}, it suffices to prove the
theorem for $\ep=1$ and $\bv(x,t)=\bv^1(x,t)$.

Since $\ell-d\le 0$, it follows from~\eqref{eqFk4} that
\begin{equation*}
\left|\dfrac{1}{t^{(k+d)/\ell-1}}
{F}_k\left(\dfrac{x}{t^{1/\ell}}\right)\right|\to 0\ \text{as }
t\to\infty,\quad k=0,1,2,\dots.
\end{equation*}
Therefore, by Lemma~\ref{lAsymphMultiDimPrelim},
$$
\Omega(x)=\lim\limits_{t\to \infty}\bv(x,t),\quad x\in\bbZ^d.
$$

On the other hand, representation~\eqref{eqy_nFourierMultiDim}
with $\ep=1$, estimate~\eqref{eqDiscreteEllipticitym},
  and Lebesgue's dominated convergence theorem imply that
$$
\lim\limits_{t\to
\infty}\bv(x,t)=\dfrac{1}{(2\pi)^d}\int_{R_\pi}\dfrac{e^{i
y\theta}}{A(\theta)}\,d\theta,\quad x\in\bbZ^d,
$$
where the integral is absolutely convergent for $\ell-d\le -1$ due to~\eqref{eqDiscreteEllipticitym}. The
last two equalities, combined with
Lemma~\ref{lAsymphMultiDimPrelim} and Remark~\ref{remcFkEqual0},
complete the proof.
\endproof

As in Corollaries~\ref{corApproxt} and~\ref{corApproxEpsilon},
Theorem~\ref{thAsympfGeneral1} and Lemma~\ref{lFk12},
item~\ref{lFk3}, imply that the higher approximation order $M$
(see Definition~\ref{defAppriximation}) we choose, the faster the
second Green function converges to the main term in its asymptotic
representation, as $t\to\infty$ or as $\ep\to0$.

\subsection{Case $\ell-d\ge 0$}\label{subsecEll-dPositive}

First, we take $\ep=1$ and fix an arbitrary $x\in\bbZ^d$. Using
representation~\eqref{eqy_nFourierMultiDim} with $\ep=1$, we write
\begin{equation}\label{eqldBigger0v1}
\bv(x,t)=\dfrac{1}{(2\pi)^d}\big(\bv_1(x,t)+\bv_2(x,t)-\bv_3(x,t)\big),
\end{equation}
\begin{equation}\label{eqldBigger0v2}
\begin{aligned}
& \bv_1(x,t):=\int_{B_{t^{-1/\ell}}} \dfrac{1-e^{-t A(\theta)}}{
A(\theta)}e^{i x\theta}\,d\theta,\quad
\bv_2(x,t):=\int_{R_\pi\setminus B_{t^{-1/\ell}}} \dfrac{e^{i
x\theta}}{ A(\theta)}\,d\theta, \\
&\bv_3(x,t):=\int_{R_\pi\setminus B_{t^{-1/\ell}}} \dfrac{e^{-t
A(\theta)}}{ A(\theta)}e^{i x\theta}\,d\theta,
\end{aligned}
\end{equation}
where
$
B_\lambda:=\{\theta\in\bbR^d: |\theta|<\lambda\}
$
and $R_\pi$ is the cube defined by~\eqref{eqRepsilon}.

{\bf 1. Estimate of $\bv_1(x,t)$.} Making the change of variables
$\theta=t^{-1/\ell}\xi$, we obtain
\begin{equation}\label{eqv1NewVar}
\bv_1(t)=t^{-d/\ell}\int_{B_1} \dfrac{1-e^{-t A(t^{-1/\ell}\xi)}}{
 A(t^{-1/\ell}\xi)} e^{i x t^{-1/\ell}\xi}\,d\xi.
\end{equation}
Using the Taylor expansion for $e^{ix\theta}$ and
formula~\eqref{eqConnectionSymbols}, we represent
\begin{equation}\label{eqldBigger0v2_1}
\dfrac{e^{ix\theta}}{A(\theta)}=\sum\limits_{k=0}^{\ell-d}
{Q}_{k-\ell}(\theta,x)+\hat{Q}(\theta,x),
\end{equation}
where ${Q}_{-\ell}(\theta,x)=\dfrac{1}{\cA(\theta)}$,
${Q}_{k-\ell}(\theta,x)$ are homogeneous functions of degree
$k-\ell$ with respect to $\theta$ representable as ratios of homogeneous polynomials, and there is $c_1=c_1(x)>0$ such
that
\begin{equation}\label{eqldBigger0v2_2}
|\hat{Q}(\theta,x)|\le c_1|\theta|^{-d+1},\quad \theta\in R_\pi.
\end{equation}
Relations~\eqref{eqldBigger0v2_1} and~\eqref{eqldBigger0v2_2}
yield
\begin{equation}\label{eqldBigger0v2_1xi}
\dfrac{e^{i x
t^{-1/\ell}\xi}}{A(t^{-1/\ell}\xi)}=\sum\limits_{k=0}^{\ell-d}
   t^{-(k-\ell)/\ell} {Q}_{k-\ell}(\xi,x)+\hat{Q}(t^{-1/\ell}\xi,x),
\end{equation}
where
\begin{equation}\label{eqldBigger0v2_2xi}
|\hat{Q}(t^{-1/\ell}\xi,x)|\le c_1 t^{(d-1)/\ell}
|\xi|^{-d+1},\quad \xi\in R_{\pi t^{1/\ell}},\ t\ge t_0.
\end{equation}

On the other hand, due to~\eqref{eqAsympAthetaMultiDim}
and~\eqref{eqAsymphWeak4'5MultiDim},
\begin{equation}\label{eqldBigger0v2_1ExpXxi}
1-e^{-t
A(t^{-1/\ell}\xi)}=1-e^{-\cA(\xi)}-e^{-\cA(\xi)}\left(\sum\limits_{m=M}^K
t^{-m/\ell}\cP_{m+\ell}(\xi)+\hat P(\xi,t)\right),
\end{equation}
 where
\begin{equation}\label{eqAsymphWeak4'5MultiDim'SmallXi}
|\hat P(\xi,t)|\le  c_2  t^{-(K+1)/\ell} |\xi|^{K+1+\ell}, \quad
\xi\in B_1,\ t\ge t_0.
\end{equation}

Combining~\eqref{eqv1NewVar} with
\eqref{eqldBigger0v2_1xi}--\eqref{eqAsymphWeak4'5MultiDim'SmallXi}
and recalling that $K\ge\ell-d$, we obtain
$$
\begin{aligned}
\bv_1(x,t) & = \sum\limits_{k=0}^{\ell-d} t^{1-(k+d)/\ell}
 \int\limits_{B_1}\left(1-e^{-\cA(\xi)}\right){Q}_{k-\ell}(\xi,x)\,d\xi
\\
& - \sum\limits_{k=0}^{\ell-d-M} \sum\limits_{m=M}^{K}
t^{1-(k+m+d)/\ell}
 \int\limits_{B_1}
 e^{-\cA(\xi)}{Q}_{k-\ell}(\xi,x)\cP_{m+\ell}(\xi)\,d\xi +
 o(1)\quad\text{as } t\to\infty.
\end{aligned}
$$
Therefore, explicitly writing a coefficient at the zeroth power of
$t$, we have
\begin{equation}\label{eqv1Final}
\begin{aligned}
\bv_1(x,t) & = \bv_1^*(x,t)+
 \int\limits_{B_1}\left(1-e^{-\cA(\xi)}\right){Q}_{-d}(\xi,x)\,d\xi
\\
& - \sum\limits_{k=0}^{\ell-d-M}
 \int\limits_{B_1}
 e^{-\cA(\xi)}{Q}_{k-\ell}(\xi,x)\cP_{2\ell-d-k}(\xi)\,d\xi +
 o(1)\quad\text{as } t\to\infty,
\end{aligned}
\end{equation}
where $\bv_1^*(x,t)$ is a linear combination of the positive
powers $t^{1-(k+d)/\ell}$ ($k=0,\dots,\ell-d-1$) with
$x$-dependent coefficients.

{\bf 2. Estimate of $\bv_2(x,t)$.} We   write
${Q}_{k-\ell}(\theta,x)$ (defined in~\eqref{eqldBigger0v2_1})
 in the spherical coordinates $(r,\ph)$ with respect to $\theta$:
\begin{equation}\label{eqldBigger0v2_3}
{Q}_{k-\ell}(\theta,x)=r^{k-\ell}\tilde{Q}_{k-\ell}(\ph,x)
\end{equation}
and note that $\tilde{Q}_{k-\ell}(\ph,x)$ are infinitely
differentiable with respect to $\ph$ and $x$. Recall the function $
r_\pi(\ph)$ defined in~\eqref{eqr_pi}.
 Then, using~\eqref{eqldBigger0v2}, \eqref{eqldBigger0v2_1}, and~\eqref{eqldBigger0v2_3},
 we obtain
\begin{equation*}
\begin{aligned}
\bv_2(x,t) &=
  \sum\limits_{k=0}^{\ell-d} \int\limits_\Phi
\tilde{Q}_{k-\ell}(\ph,x){J}(\ph)\,d\ph\int\limits_{t^{-1/\ell}}^{r_\pi(\ph)}r^{d+k-\ell-1}\,dr + \int\limits_{R_\pi\setminus
B_{t^{-1/\ell}}} \hat{Q}(\theta,x)\, d\theta.
\end{aligned}
\end{equation*}
Since $\hat{Q}(\theta,x)$ is integrable at the origin due
to~\eqref{eqldBigger0v2_2}, we calculate (as $t\to\infty$)
\begin{equation}\label{eqv2Final}
\begin{aligned}
\bv_2(x,t) &=  \bv_2^*(x,t)  + \dfrac{\ln
t}{\ell}\int\limits_\Phi \tilde{Q}_{-d}(\ph,x){J}(\ph)\,d\ph \\
& - \sum\limits_{k=0}^{\ell-d-1}
\dfrac{1}{\ell-d-k}\int\limits_\Phi (r_\pi(\ph))^{d+k-\ell}\tilde{Q}_{k-\ell}(\ph,x){J}(\ph)\,d\ph +
\int\limits_{R_\pi} \hat{Q}(\theta,x)\, d\theta \\
& +  \int\limits_\Phi \ln(r_\pi(\ph))\tilde{Q}_{-d}(\ph,x){J}(\ph)\,d\ph + o(1),
\end{aligned}
\end{equation}
where $\bv_2^*(x,t)$ is a linear combination of the positive
powers $t^{1-(k+d)/\ell}$ ($k=0,\dots,\ell-d-1$) with
$x$-dependent coefficients.

{\bf 3. Estimate of $\bv_3(x,t)$.} Making the change of variables
$\theta=t^{-1/\ell}\xi$, we obtain from~\eqref{eqldBigger0v2}
\begin{equation*}
\bv_3(t)=t^{-d/\ell}\int_{R_{\pi t^{1/L}}\setminus B_1}
\dfrac{e^{-t A(t^{-1/\ell}\xi)}}{
 A(t^{-1/\ell}\xi)} e^{i x t^{-1/\ell}\xi}\,d\xi + t^{-d/\ell}\int_{R_{\pi t^{1/\ell}}\setminus R_{\pi
t^{1/L}}} \dfrac{e^{-t A(t^{-1/\ell}\xi)}}{
 A(t^{-1/\ell}\xi)} e^{i x t^{-1/\ell}\xi}\,d\xi,
\end{equation*}
where $L$ is defined by~\eqref{eqAsymphWeak4'MultidimL}.
Therefore, using~\eqref{eqDiscreteEllipticitym}, we obtain
similarly to~\eqref{eqAsymphWeak7MultiDim}
\begin{equation}\label{eqv3NewVar}
\bv_3(t)=t^{-d/\ell}\int_{R_{\pi t^{1/L}}\setminus B_1}
\dfrac{e^{-t A(t^{-1/\ell}\xi)}}{
 A(t^{-1/\ell}\xi)} e^{i x t^{-1/\ell}\xi}\,d\xi +
 o(1)\quad\text{as } t\to\infty.
\end{equation}

Using~\eqref{eqv3NewVar}, \eqref{eqAsympAthetaMultiDim},
\eqref{eqAsymphWeak4'5MultiDim},  and \eqref{eqldBigger0v2_1xi},
we have
$$
\begin{aligned}
\bv_3(x,t) & =t^{-d/\ell}\int_{R_{\pi t^{1/L}}\setminus B_1}
e^{-\cA(\xi)}\left(1+ \sum\limits_{m=M}^K
t^{-m/\ell}\cP_{m+\ell}(\xi)+
 \hat P(\xi,t)\right)\\
 & \cdot\left(\sum\limits_{k=0}^{\ell-d}
   t^{-(k-\ell)/\ell}
   {Q}_{k-\ell}(\xi,x)+\hat{Q}(t^{-1/\ell}\xi,x)\right)\,d\xi +
   o(1).
 \end{aligned}
$$
Taking into account estimates~\eqref{eqAsymphWeak4'5MultiDim'}
and~\eqref{eqldBigger0v2_2xi} and recalling that $K\ge\ell-d$, we obtain
\begin{equation}\label{eqv3Final}
\begin{aligned}
&\bv_3(x,t)  = \int_{\bbR^d\setminus B_1} e^{-\cA(\xi)}\left(1+
\sum\limits_{m=M}^K t^{-m/\ell}\cP_{m+\ell}(\xi)\right)
\left(\sum\limits_{k=0}^{\ell-d}
   t^{1-(k+d)/\ell}
   {Q}_{k-\ell}(\xi,x)\right)\,d\xi +
   o(1)\\
 &  =\bv_3^*(x,t)+ \int_{\bbR^d\setminus
 B_1}e^{-\cA(\xi)}{Q}_{-d}(\xi,x)\,d\xi + \sum\limits_{k=0}^{\ell-d-M}
 \int\limits_{\bbR^d\setminus B_1}
 e^{-\cA(\xi)}{Q}_{k-\ell}(\xi,x)\cP_{2\ell-d-k}(\xi)\,d\xi+o(1),
 \end{aligned}
\end{equation}
where $\bv_3^*(x,t)$ is a linear combination of the positive
powers $t^{1-(k+d)/\ell}$ ($k=0,\dots,\ell-d-1$) with
$x$-dependent coefficients.

{\bf 4. Estimate of $\bv(x,t)$.} Combining~\eqref{eqldBigger0v1}
with \eqref{eqv1Final}, \eqref{eqv2Final}, \eqref{eqv3Final}, we
have for each $x\in\bbZ^d$
\begin{equation}\label{eqldBigger0v11}
\begin{aligned}
\bv(x,t) & = \bv^*(x,t)+ \dfrac{\ln
t}{(2\pi)^d\ell}\int\limits_\Phi
\tilde{Q}_{-d}(\ph,x){J}(\ph)\,d\ph + {S}(x)  + o(1),
\end{aligned}
\end{equation}
where $\bv^*(x,t)$ is a linear combination of the positive powers
$t^{1-(k+d)/\ell}$ ($k=0,\dots,\ell-d-1$) with $x$-dependent
coefficients and
\begin{equation}\label{eqldBigger0v12}
\begin{aligned}
{S}(y)  :=\dfrac{1}{(2\pi)^d}  \Bigg( &
 \int\limits_{B_1}\left(1-e^{-\cA(\xi)}\right){Q}_{-d}(\xi,y)\,d\xi -
\int_{\bbR^d\setminus
 B_1}e^{-\cA(\xi)}{Q}_{-d}(\xi,y)\,d\xi
\\
& + \int\limits_{R_\pi} \hat{Q}(\theta,y)\, d\theta -
\sum\limits_{k=0}^{\ell-d-M}
 \int\limits_{\bbR^d}
 e^{-\cA(\xi)}{Q}_{k-\ell}(\xi,y)\cP_{2\ell-d-k}(\xi)\,d\xi   \\
& - \sum\limits_{k=0}^{\ell-d-1}
\dfrac{1}{\ell-d-k}\int\limits_\Phi (r_\pi(\ph))^{d+k-\ell}\tilde{Q}_{k-\ell}(\ph,y){J}(\ph)\,d\ph  \\
& +  \int\limits_\Phi \ln(r_\pi(\ph))\tilde{Q}_{-d}(\ph,y){J}(\ph)\,d\ph   \Bigg).
\end{aligned}
\end{equation}

Now we are in a position to formulate our main result in the case
$\ell-d\ge 0$.

\begin{theorem}\label{thAsympfGeneral2}
Let $\ell-d\ge 0$. Then, for any $\ep>0$, $t_0>0$, integer $
K\ge\max(M,\ell-d)$, and all $t\ge t_0\ep^\ell$, the following
hold.
\begin{enumerate}
\item If $x\in \bbR_\ep^d\setminus\{0\}$, then
\begin{equation}\label{eqAsympfGeneralMultiDim-d-lBigger1Final1}
\begin{aligned}
\bv^\ep(x,t) & = t^{1-d/\ell}\,
{F}_0\left(\dfrac{x}{t^{1/\ell}}\right) +
\sum\limits_{k=M}^{\ell-d} \ep^k t^{1-(k+d)/\ell}\,
{F}_k\left(\dfrac{x}{t^{1/\ell}}\right) + \ep^{\ell-d}\,
\Omega\left(\dfrac{x}{\ep}\right) \\
& + \sum\limits_{k=\max(M,\ell-d+1)}^K \dfrac{\ep^k}{
t^{(k+d)/\ell-1}}\, {F}_k\left(\dfrac{x}{t^{1/\ell}}\right) +
\br_\bv^\ep(K; x,t),
\end{aligned}
\end{equation}
where ${F}_k(y)$ are given by~\eqref{eqcFLog-Corrected1}
and~\eqref{eqcFLog-Corrected2}, respectively,
$$
\Omega(y)={S}(y)-\hat{F}(y) -
\sum\limits_{k=0}^{\ell-d}\Omega_k(y),
$$
${S}(y)$ is given by~\eqref{eqldBigger0v12}, $\hat{F}(y)$
by~\eqref{eqFindGamma5'}, and $ \Omega_k(y)$
by~\eqref{eqFindGamma4'}.

\item If $x=0$, then
\begin{equation}\label{eqFindGamma5X0Final}
\begin{aligned}
\bv^\ep(0,t) & = \sum\limits_{k=0}^{\ell-d-1}  \ep^k t^{1-(d+k)/\ell} \dfrac{\ell\,
{H}_k(0)}{\ell-d-k} + \ep^{\ell-d}\big({H}_{\ell-d}(0)\ln t - \ell{H}_{\ell-d}(0)\ln \ep  + {S}(0)\big)\\
&  + \sum\limits_{k=\ell-d+1}^K \dfrac{\ep^k }{t^{(d+k)/\ell-1}}
\dfrac{\ell\, {H}_k(0)}{\ell-d-k} + \br_\bv^\ep(K; 0,t),
\end{aligned}
\end{equation}
where ${S}(0)$ is given by~\eqref{eqldBigger0v12} with $y=0$.
\end{enumerate}
In both cases,
$$
|\br_\bv^\ep(K; x,t)|\le
\dfrac{\ep^{K+1}R_\bv(K,t_0)}{t^{(K+d+1)/\ell-1}}
$$
and $R_\bv(K,t_0)\ge 0$ does not depend on $\ep>0$,
$x\in\bbR_\ep^d$, and $t\ge t_0 \ep^\ell$.
\end{theorem}
\proof The result follows from Lemmas~\ref{lAsymphMultiDimPrelim}
and \ref{lAsymphMultiDimPrelimVer2},
Corollary~\ref{corAsymphMultiDimPrelimVer2}, and
formula~\eqref{eqldBigger0v11}.
\endproof

\begin{remark}\label{remAsympDerivatives}
Using Theorems~\ref{thAsymphMultiDim} and~\ref{thAsympfGeneral2},
one can easily obtain asymptotic formulas for the first-order
``spatial derivatives'' $\ep^{-1}\delta_{\ep,k\pm}\bu^\ep(x,t)$
and $\ep^{-1}\delta_{\ep,k\pm}\bv^\ep(x,t)$, $k=1,\dots,d$, and
for their higher-order analogues. To do so, one should use
Taylor's expansions around $x/t^{1/\ell}$ of the functions
entering formulas~\eqref{eqAsymphGeneralMultiDim}
and~\eqref{eqAsympfGeneralMultiDim-d-lBigger1Final1}. We will not
provide this result in full generality, but rather present it for
a one-dimensional example in
Section~\ref{subsec1DExampleGradient}.
\end{remark}

\section{Higher order approximations of the 1D Laplacian}\label{sec1DExample}

In this section, we illustrate our general results by applying them to the
one-dimensional discrete Laplace operator.
Throughout this section, we fix  $N\in\bbN$. In what
follows, we will not explicitly indicate the dependence of
constants, functions, etc. on $N$.

\subsection{List of constants and functions}

For asymptotic formulas below, we will explicitly find all the
coefficients. In this subsection, we collect the constants and the
functions which are needed for that.

 For integer $J\ge 0$, set
\begin{align}
& b_n :=\dfrac{2(-1)^{n+1}}{(2(n+1))!}\sum\limits_{\nu=1}^N a_\nu
\nu^{2(n+1)},\quad & & n\ge N, \label{eqbk} \\
& c_{nm} :=\sum\limits_{\substack{l_1,\dots,l_m\ge N,\\
l_1+\dots+l_m=n}}
 b_{l_1}\cdot\ldots\cdot b_{l_m} , & & n\ge N,\ m\ge 1,
\label{eqck} \\
& \cP_n(\xi):=\sum\limits_{1\le m\le n/N}\dfrac{c_{n
m}\xi^{2(n+m)}}{m!}, & & n\ge
N,\label{eqPkxi} \\
& d_{J n}:=\sum\limits_{1\le m\le \min(n/N,J)}(-1)^m \dbinom{J}{m}
c_{n m} , &
& n\ge N,\label{eqdk} \\
& \ccQ_{J n}(\xi):=0, & & n=N,\dots,2N-1,\label{eqQkxi12N}\\
& \ccQ_{J n}(\xi):=\sum\limits_{\substack{l_1,l_2\ge N,\\
l_1+l_2=n}}
d_{Jl_1}\xi^{2 l_1}\cP_{l_2}(\xi), & & n \ge 2N,\label{eqQkxiBigger2N}\\
& \cR_{J n}(\xi):=(-1)^J\xi^{2J}(d_{J
n}\xi^{2n}+\cP_n(\xi)+\ccQ_{J n}(\xi)), & & n\ge N.\label{eqRkxi}
\end{align}
As before,  we use the convention  that the sum vanishes if the
set over which the summation is taken is empty.

\subsection{Setting}

Let $\bu^\ep(x)$, $x\in\bbR_\ep^1$, be a grid function. We consider the
difference operator
\begin{equation}\label{eqDeltaEp}
\bA_\ep=-\bDelta_\ep,\quad \bDelta_\ep
\bu(x):=\ep^{-2}\sum\limits_{\nu=1}^N
a_\nu\big(\bu(x-\ep\nu)-2\bu(x)+\bu(x+\ep\nu)\big),
\end{equation}
with the coefficients $a_\nu$ satisfying the equalities
\begin{equation}\label{eqanu}
\sum\limits_{\nu=1}^N a_\nu \nu^2=1,\qquad  \sum\limits_{\nu=1}^N
a_\nu \nu^{2m}=0,\quad m=2,\dots,N.
\end{equation}
The matrix of system~\eqref{eqanu} is the Vandermonde matrix,
which guarantees the existence and uniqueness of $a_\nu$. Thus,
throughout this section, $\ell=2$ and $d=1$.


The following lemma shows that Condition~\ref{condAppriximation}
holds with $M=2N$.
\begin{lemma}\label{lem1DApprox}
The operator $\bDelta_\ep$ is the $(2N)$th order approximation of
the second derivative $($the one-dimensional Laplacian$)$ in the
sense of Definition~$\ref{defAppriximation}$.
\end{lemma}
\proof For any smooth function $u(y)$ and $\ep>0$,
equalities~\eqref{eqanu} imply
\begin{equation*}
\begin{aligned}
& \ep^{-2}\sum\limits_{\nu=1}^N
a_\nu\big(u(y-\ep\nu)-2u(y)+u(y+\ep\nu)\big)\\
& = \ep^{-2} \sum\limits_{\nu=1}^N a_\nu\left(\sum\limits_{m=1}^N
\dfrac{2(\ep\nu)^{2m}}{(2m)!}
u^{(2m)}(y)+\dfrac{(\ep\nu)^{2N+2}}{(2N+2)!}\left(u^{(2N+2)}(y-\xi_{1\nu}\ep\nu)+u^{(2N+2)}(y+\xi_{2\nu}\ep\nu)\right)\right)\\
&=\ep^{-2}\left(\sum\limits_{m=1}^N
u^{(2m)}(y)\dfrac{2\ep^{2m}}{(2m)!} \sum\limits_{\nu=1}^N a_\nu
\nu^{2m}\right)  +\ep^{2N} R_N(y) =u''(y)+\ep^{2N}R_N(y),
\end{aligned}
\end{equation*}
where
\begin{equation*}
R_N(y)=\sum\limits_{\nu=1}^N
\dfrac{\nu^{2N+2}}{(2N+2)!}\left(u^{(2N+2)}(y-\xi_{1\nu}\ep\nu)+u^{(2N+2)}(y+\xi_{2\nu}\ep\nu)\right),
\quad \xi_{1\nu},\xi_{2\nu}\in[0,1].
\end{equation*}
\endproof

Due to Definition~\ref{defSymbolDiscrete},  the symbol of the
operator $\bA_\ep=-\bDelta_\ep$ is given by
\begin{equation}\label{eqALaplace2N}
 A(\theta)=2\sum\limits_{\nu=1}^N
a_\nu(1-\cos\nu\theta).
\end{equation}

The following lemma shows that
Condition~\ref{condDiscreteEllipticity} holds. The proof is given
in Appendix~\ref{appendProofEllipt}.

\begin{lemma}\label{lAtheta}
Let $A(\theta)$ be given by~\eqref{eqALaplace2N}. Then
\begin{enumerate}
\item $A'(0)=A'(\pi)=0$ and $A'(\theta)>0$ for all
$\theta\in(0,\pi)$$;$

\item  $A(0)=0$ and $A(\theta)>0$ for all $\theta\in (0,\pi]$$;$

\item\label{lAtheta3} $A(\theta)=\theta^2- \sum\limits_{n=N}^K
{b_n}\theta^{2n+2}+O(\theta^{2(K+2)}), $ where $b_n$ are given
by~\eqref{eqbk}.
\end{enumerate}
\end{lemma}

From now on, we confine ourselves by the case $\ep=1$. The case
$\ep>0$ can be easily obtained by applying the scaling
rule~\eqref{eqRescalingGreen}. Let $\bdelta(x):=\bdelta^1(x)$,
$x\in\bbZ$, be the grid delta-function defined
in~\eqref{eqdeltaEp} with $\ep=1$.

\begin{definition}\label{defGreen1D}
Let $\ep=1$. We call the function $\bu(x,t)$, $x\in\bbZ,\ t\ge 0$,
the {\em first discrete Green function} if $\bu(\cdot,t)$ is a
rapidly decreasing grid function for all $t\ge0$, $\bu(x,\cdot)\in
C^1[0,\infty)$ for all $x\in\bbZ$, and
$$
\left\{
\begin{aligned}
& \dot \bu(x,t)-\bDelta_1 \bu(x,t)=0, & & x\in\bbZ,\ t>0,\\
& \bu(x,0)=\bdelta(x),   & & x\in\bbZ.
\end{aligned}
\right.
$$

We call the function $\bv(x,t)$, $x\in\bbZ,\ t\ge 0$, the {\em
second discrete Green function} if $\bv(\cdot,t)$ is a rapidly
decreasing grid function for all $t\ge0$, $\bv(x,\cdot)\in
C^1[0,\infty)$ for all $x\in\bbZ$, and
$$
\left\{
\begin{aligned}
& \dot \bv(x,t)-\bDelta_1 \bv(x,t)=\bdelta(x), & & x\in\bbZ,\ t>0,\\
& \bv(x,0)=0,   & & x\in\bbZ.
\end{aligned}
\right.
$$
\end{definition}


\subsection{Asymptotics of the first Green function
$\bu(x,t)$}

We begin with asymptotic formulas for the first Green function
$\bu(x,t)$. Consider the function
\begin{equation}\label{eqhtildeh}
 h(y):=\dfrac{1}{2\sqrt{\pi}}\,e^{-y^2/4},
\end{equation}
which coincides with ${H}(y)$ defined in~\eqref{eqhGeneral}. Note
that $\dfrac{1}{\sqrt{t}}h(y)$ is the Green function of the
operator $u_t(y,t)-u_{yy}(y,t)$ (cf.~\eqref{eqContGreen}).

\begin{theorem}\label{thAsymph}
For any $t_0>0$, integer $J\ge 0$ and ${K_1} \ge N$, and all
$x\in\bbZ$ and $t\ge t_0$, we have
\begin{equation}\label{eqAsymphGeneral}
\dfrac{d^{J}\bu(x,t)}{dt^J}=\dfrac{1}{t^J\sqrt{t}}
h_{J0}\left(\dfrac{x}{\sqrt{t}}\right)+\sum\limits_{n=N}^ {K_1}
\dfrac{1}{t^{n+J}\sqrt{t}}\, h_{J
n}\left(\dfrac{x}{\sqrt{t}}\right) + \br_\bu(J, {K_1};x,t),
\end{equation}
where
\begin{equation}\label{eqhJk}
h_{J0}(x):=\dfrac{d^{2J}h(x)}{dx^{2J}},\quad h_{J n}(x):=\cR_{J
n}\left(-i\dfrac{d}{dx}\right) h(x),
\end{equation}
$$
|\br_\bu(J, {K_1};x,t)|\le \dfrac{R_\bu(J, {K_1} ,t_0)}{t^{J+
{K_1} +1}\sqrt{t}},
$$
$\cR_{J n}(\xi)$ are polynomials given by~\eqref{eqRkxi}, and $R_\bu(J, {K_1}
,t_0)\ge 0$ does not depend on $t\ge t_0$ and $x\in\bbZ$.
\end{theorem}
\proof Due to Lemma~\ref{lAtheta}, part~\ref{lAtheta3}, the Taylor
expansion of the symbol $A(\theta)$ contains only even powers of
$\theta$. This fact and Theorem~\ref{thAsymphMultiDim} with
$M=2N$, $K=2 K_1$, and $k=2n$ imply the result.
\endproof

\subsection{Asymptotics of the second Green function
$\bv(x,t)$}

  For
$y\ge 0$, we introduce the functions
\begin{equation}\label{eqf}
f_0(y)=2y\int_y^\infty \rho^{-2}h(\rho)\,d\rho,
\end{equation}
\begin{equation}\label{eqFormulafk}
f_n(y)=-\dfrac{2}{y^{2n-1}}\int_0^y
\rho^{2n-2}h_{0n}(\rho)\,d\rho,\quad n=1,2,\dots.
\end{equation}

Note that
\begin{equation}\label{eqf''h}
f_0''(y)=h(y),\quad y>0,
\end{equation}
where $h(y)$ is given by~\eqref{eqhtildeh}. Furthermore, $f_0(y)$
is positive, real analytic for\footnote{Analyticity at $y=0$ is
understood in the sense that $f(y)$ can be represented as the
Taylor series with respect to the powers of $y$ which converges to
$f(y)$ in a right-hand neighborhood of $y=0$.} $y\ge 0$ (e.g.,
use~\eqref{eqf''h}), and vanishes at infinity, together with all
its derivatives.


Using~\eqref{eqhtildeh} and~\eqref{eqhJk}, we see that $f_n(y)$
are  real analytic for  $y\ge 0$ and vanish at infinity, together
with all their derivatives.

\begin{theorem}\label{thAsympfGeneral}
For any $t_0>0$, integer $  {K_1} \ge N$, and all $x\in\bbZ$ and
$t\ge t_0$, we have
\begin{equation}\label{eqAsympfGeneral}
\begin{aligned}
 \bv(x,t)&=\sqrt{t}\, f_0\left(\dfrac{x}{\sqrt{t}}\right)+\Omega(x)+\sum\limits_{n=N}^ {K_1}  \dfrac{1}{t^{n-1}\sqrt{t}}\,
 f_n\left(\dfrac{x}{\sqrt{t}}\right)+
\br_\bv( {K_1};x,t),
& & x\ge0,\\
\bv(x,t)&=\bv(-x,t), & & x\le -1,
\end{aligned}
\end{equation}
where
\begin{equation}\label{eqGamma1D}
\Omega(x)=\dfrac{1}{2\pi}\left(\,\int\limits_{-\pi}^\pi
\left(\dfrac{\cos(x\xi)}{A(\xi)}-\dfrac{1}{\xi^2}\right)d\xi -
\dfrac{2}{\pi} \right) + \dfrac{x}{2},
\end{equation}
$f_0(y), f_1(y),\dots,f_ {K_1} (y)$ are given by~\eqref{eqf},
\eqref{eqFormulafk},
$$
|\br_\bv( {K_1};x,t)|\le \dfrac{R_\bv( {K_1} ,t_0)}{t^{ {K_1}
}\sqrt{t}},
$$
and $R_\bv( {K_1} ,t_0)\ge 0$ does not depend on $t\ge t_0$ and
$x\ge0$.

Moreover, if $N=1$, then $\Omega(x)=0$ for all integer $x\ge 0$.
\end{theorem}
\proof Since $v(x,t)\equiv v(-x,t)$, $x\in\bbZ$, it remains to
prove the first formula in~\eqref{eqAsympfGeneral}. Due to
Lemma~\ref{lAtheta}, part~\ref{lAtheta3}, the Taylor expansion of
the symbol $A(\theta)$ contains only even powers of $\theta$.
Therefore, ${H}_k(y)\equiv{H}_{0k}(y)\equiv 0$ for odd $k$. Hence,
${F}_k(y)\equiv 0$ for odd $k$, where ${F}_k(y)$ are given
by~\eqref{eqcFLog-Corrected2}. On the other hand,
relations~\eqref{eqf}, \eqref{eqFormulafk}
and~\eqref{eqcFLog-Corrected1}, \eqref{eqcFLog-Corrected2}, imply
that
\begin{equation}\label{eqfEqualscF}
f_0(y)\equiv{F}_0(y),\qquad f_n(y)\equiv{F}_{2n}(y),\quad n\ge 1.
\end{equation}

Thus, applying Theorem~\ref{thAsympfGeneral2}   with $\ep=1$,
$M=2N$,   $K=2 K_1$, and $k=2n$, we obtain the first formula
in~\eqref{eqAsympfGeneral}
 with $\Omega(x)$ given by
\begin{equation}\label{eqGamma1D_1}
\Omega(x)=\dfrac{1}{2\pi}\left(\,\int\limits_{-\pi}^\pi
\left(\dfrac{\cos(x\xi)}{A(\xi)}-\dfrac{1}{\xi^2}\right)d\xi -
\dfrac{2}{\pi} \right) - 2x \int\limits_0^1
\dfrac{h(\rho)-h(0)}{\rho^2}\,d\rho - 2x\left(\int\limits_1^\infty
\dfrac{h(\rho)}{\rho^2}\,d\rho-h(0)\right).
\end{equation}
Note that
$$
- 2 \int\limits_0^1 \dfrac{h(\rho)-h(0)}{\rho^2}\,d\rho -
2\left(\int\limits_1^\infty
\dfrac{h(\rho)}{\rho^2}\,d\rho-h(0)\right)=-\dfrac{1}{\sqrt{\pi}}\int\limits_0^\infty
\dfrac{e^{-\rho^2/4}-1}{\rho^2}\,d\rho=\dfrac{1}{2},
$$
where the last equality can be proved, e.g., by introducing a
parameter $a$ as follows: $e^{-a\rho^2/4}$ and differentiating
with respect to it. Hence,~\eqref{eqGamma1D_1}
yields~\eqref{eqGamma1D}.

 Finally, if $N=1$, then $\Omega(x)=0$ for all $x\in\bbZ$ by
Lemma~\ref{lGammaIs0}.
\endproof

\subsection{Asymptotics of the gradient of the Green  functions}\label{subsec1DExampleGradient}

Next, we derive asymptotic formulas for $\nabla
\bv(x,t):=\bv(x+1,t)-\bv(x,t)$ and its time derivatives, see
Remark~\ref{remAsympDerivatives}. We introduce the function
\begin{equation}\label{eqgtildeg}
 g(x):=f_0'(x),\quad x\ge 0,
\end{equation}
where $f_0(x)$ is given by~\eqref{eqf}. Note that $g(x)$ is
negative, real analytic for $x\ge 0$,  and vanishes at infinity,
together with all its derivatives.

Furthermore, for $n\ge 2$, we set
\begin{equation}\label{eqgk}
 g_n(x):= \dfrac{1}{(n+1)!}f_0^{(n+1)}(x)+\sum\limits_{\substack{s\ge N,\,j\ge 1,\\
 2s+j-1=n}}\dfrac{1}{j!}f_s^{(j)}(x),\quad x\ge 0,
\end{equation}
where $f_s(x)$ are given by~\eqref{eqFormulafk}. Note that
$f_s(x)$ are real analytic for $x\ge 0$  and vanish at infinity,
together with all their derivatives. Therefore, expanding
 the functions $f_s$ around
$x/\sqrt{t}$ in~\eqref{eqAsympfGeneral}, we arrive at the
following theorem.

\begin{theorem}\label{thAsympgGeneral}
For any $t_0>0$, integer $ {K_1} \ge 2$, and all $x\in\bbZ$ and
$t\ge t_0$, we have
\begin{equation}\label{eqAsympgGeneral}
\begin{aligned}
\nabla \bv(x,t)&=\nabla \Omega(x)+
g\left(\dfrac{x}{\sqrt{t}}\right) +\dfrac{1}{2\sqrt{t}}\,
h\left(\dfrac{x}{\sqrt{t}}\right) & &\\
&+\sum\limits_{n=2}^ {K_1} \dfrac{1}{t^{n/2}}\,
g_n\left(\dfrac{x}{\sqrt{t}}\right)+ \hat\br_\bv( {K_1};x,t),
& & x\ge0,\\
 \nabla \bv(x,t)&=-\nabla \bv(-(x+1),t), & & x\le -1.
\end{aligned}
\end{equation}
where $\Omega(x)$ is given by~\eqref{eqGamma1D}, $g_2(x),\dots,g_
{K_1} (x)$ are given by~\eqref{eqgk},
$$
|\hat\br_\bv( {K_1};x,t)|\le \dfrac{\hat R_\bv( {K_1}
,t_0)}{t^{( {K_1} +1)/2}},
$$
and $\hat R_\bv( {K_1} ,t_0)\ge 0$ does not depend on $t\ge t_0$
and $x\ge0$.

Moreover, if $N=1$, then $\nabla\Omega(x)=0$ for all integer $x\ge
0$.
\end{theorem}

Finally,   for $n\ge 2$, we set
\begin{equation}\label{eqgJk}
 g_{J n}(x):= \dfrac{1}{(n+1)!}h_{J0}^{(n+1)}(x)+\sum\limits_{\substack{s\ge N,\,j\ge 1,\\
 2s+j-1=n}}\dfrac{1}{j!}h_{Js}^{(j)}(x),\quad x\ge 0,
\end{equation}
where $h_{J0}(y)$ and $h_{Js}(y)$ are given by~\eqref{eqhJk}. Note
that $h_{J0}(y)$ and $h_{Js}(y)$  are real analytic for $x\ge 0$
and vanish at infinity, together with all their derivatives.
Therefore, expanding
 the functions $h_{J0}$ and $h_{Js}$ around
$x/\sqrt{t}$ in~\eqref{eqAsymphGeneral}, we arrive at the
following theorem.

\begin{theorem}\label{thAsympgGeneralJ}
For any $t_0>0$, integer $J\ge 0$ and $ {K_1} \ge 2$, and all
$x\in\bbZ$ and  $t\ge t_0$, we have
\begin{equation}\label{eqAsympgGeneralJ}
\begin{aligned}
\dfrac{d^{J}\nabla\bu(x,t)}{dt^{J}} &=
\dfrac{1}{t^{J+1}}h^{(2J+1)}(x) +
\dfrac{1}{2t^{J+1}\sqrt{t}}h^{(2J+2)}(x) & & \\
& + \sum\limits_{n=2}^ {K_1} \dfrac{1}{t^{J+1+n/2}}\,
g_{J n}\left(\dfrac{x}{\sqrt{t}}\right)+ \hat\br_\bu(J, {K_1};x,t),& & x\ge 0,\\
\dfrac{d^{J}\nabla\bu(x,t)}{dt^{J}} &
=-\dfrac{d^{J}\nabla\bu(-(x+1),t)}{dt^{J}}, & & x<0.
\end{aligned}
\end{equation}
where $g_{J2}(x),\dots,g_{J{K_1}}(x)$ are given
by~\eqref{eqgJk}$;$ furthermore,
$$
|\hat\br_\bu(J, {K_1};x,t)|\le \dfrac{\hat R_\bu(J, {K_1}
,t_0)}{t^{J+1+( {K_1} +1)/2}},
$$
and $\hat R_\bu(J, {K_1} ,t_0)\ge 0$ does not depend on $t\ge t_0$
and $x\ge0$.
\end{theorem}

\appendix

\section{Proof of Lemma~$\ref{lAtheta}$ (Ellipticity)}\label{appendProofEllipt}

Let us prove assertion 1.

{\bf Step 1.} Let $\cT_\nu(y)$ and $\cU_\nu(y)$ be the Chebyshev
polynomials of the first and the second kind, respectively. Using
the relations
$$
\cos\nu\theta=\cT_\nu(\cos\theta),\quad
\dfrac{\sin\nu\theta}{\sin\theta}=\cU_{\nu-1}(\cos\theta),\quad
\nu \cU_{\nu-1}(y)=\cT_\nu(y),\quad \nu=1,2,\dots,
$$
we obtain
$$
A'(\theta)= 2 \sin\theta \sum\limits_{\nu=1}^N a_\nu
\cT_\nu'(\cos\theta).
$$
Thus, it suffices to prove that
\begin{equation}\label{eqAtheta1}
\cB(y):= \sum\limits_{\nu=1}^N a_\nu \cT_\nu'(y)>0,\quad y\in
[-1,1].
\end{equation}

{\bf Step 2.} Suppose we have proved that
\begin{equation}\label{eqAtheta2}
(-1)^{j-1}   \cB^{(j-1)}(1)>0,\quad j=1,\dots,N.
\end{equation}

Since $\cT_\nu(y)$ is a polynomial of degree $\nu$, it follows
from~\eqref{eqAtheta2} that
\begin{equation}\label{eqAtheta4}
(-1)^{N-1}\dfrac{d^{N-1}
\cB(y)}{dy^{N-1}}=(-1)^{N-1}\dfrac{d^{N-1}
\cB(1)}{dy^{N-1}}>0,\quad y\in[0,1].
\end{equation}
Therefore, $(-1)^{N-2}\dfrac{d^{N-2} \cB(y)}{dy^{N-2}}$ is
decreasing and, due to~\eqref{eqAtheta2}, is positive for
$y\in[-1,1]$. Hence, $(-1)^{N-3}\dfrac{d^{N-3} \cB(y)}{dy^{N-3}}$
is decreasing and, due to~\eqref{eqAtheta2}, is also positive for
$y\in[-1,1]$. Continuing by induction, we conclude that $\cB(y)$
is positive for $y\in[-1,1]$.

{\bf Step 3.} It remains to prove~\eqref{eqAtheta2}, i.e.,
\begin{equation}\label{eqAtheta5}
(-1)^{j-1}\sum\limits_{\nu=1}^N a_\nu \cT_\nu^{(j)}(1)>0,\quad
j=1,\dots,N.
\end{equation}

 Set
$$
u(\theta):=\cos(\theta),\quad v(y):=\arccos(y).
$$
Then
$\cT_\nu(y)=u(\nu v(y))$.
In what follows, we will use the following representations for
$\alpha=0,1,2,\dots$:
\begin{equation}\label{eqAtheta6}
v^{(\alpha)}(y)=(-1)^{\alpha}(1-y)^{-\alpha+1/2}\sum\limits_{n=0}^\infty
A_{\alpha n}(1-y)^n,
\end{equation}
where $A_{\alpha n}>0$. This series and all the series below
converge in a neighborhood of $y=1$. To prove~\eqref{eqAtheta6},
one can write $\arccos(y)=\int_y^1(1-z)^{-1/2}(1+z)^{-1/2} dz$ and
expand $(1+z)^{-1/2}$ into the Taylor series around $z=1$.

We fix some $j\in\{1,\dots,N\}$. From now on, we will not
explicitly indicate the dependence of emerging coefficients on
$j$. Since $\cT_\nu(y)=u(\nu v(y))$, we see that
\begin{equation}\label{eqAtheta6'}
\cT_\nu^{(j)}(y)= \sum\limits_{k=1}^j   \nu^k u^{(k)}(\nu
v(y))\sum\limits_{\substack{1\le l_1\le\dots\le l_k,\\
l_1+\dots+l_k=j}}B_{l_1\dots l_k}v^{(l_1)}(y)\cdot\ldots\cdot
v^{(l_k)}(y),
\end{equation}
where $B_{l_1\dots l_k}>0$. Therefore, using~\eqref{eqAtheta6}, we
have
\begin{equation}\label{eqAtheta7}
\cT_\nu^{(j)}(y)= (-1)^j \sum\limits_{k=1}^j
\sum\limits_{l=0}^\infty C_{kl}  \nu^k u^{(k)}(\nu v(y))
(1-y)^{l-j+k/2},
\end{equation}
where $C_{kl}>0$.

Expanding $u^{(k)}(\nu v(y))$ and the powers of $v(y)$
(using~\eqref{eqAtheta6} with $\alpha=0$), we obtain
\begin{equation}\label{eqAtheta8}
\begin{aligned}
& u^{(k)}(\nu v(y))=\sum\limits_{m=0}^\infty D_{km}\nu^m v^m(y),\\
& v^m(y)=(1-y)^{m/2}\sum\limits_{n=0}^\infty E_{mn}(1-y)^n,
\end{aligned}
\end{equation}
where $D_{km}\in\bbR$, $D_{km}=0$ for $k+m$ odd, and $E_{mn}>0$.
Therefore,
\begin{equation}\label{eqAtheta9}
\begin{aligned}
& u^{(k)}(\nu
v(y))=\sum\limits_{m=0}^{2j-k}\sum\limits_{n=0}^\infty
D_{km}E_{mn} \nu^m  (1-y)^{n+m/2} + (1-y)^{j-k/2} U_{k\nu}(y),
\end{aligned}
\end{equation}
where $U_{k\nu}(y)$ is continuous at $y=1$ and
\begin{equation}\label{eqAtheta10}
U_{k\nu}(1)=0.
\end{equation}

Now we substitute~\eqref{eqAtheta9} into~\eqref{eqAtheta7} and
take into account~\eqref{eqAtheta10}:
\begin{equation}\label{eqAtheta11}
\begin{aligned}
\cT_\nu^{(j)}(y)& = (-1)^j \sum\limits_{k=1}^j
\sum\limits_{m=0}^{2j-k} \sum\limits_{l,n=0}^\infty C_{kl}  D_{km}
E_{mn}   \nu^{k+m}
 (1-y)^{l+n-j+(k+m)/2}  + \cU_\nu(y),
\end{aligned}
\end{equation}
where  $U_{\nu}(y)$ is continuous at $y=1$  and
\begin{equation}\label{eqAtheta12}
U_{\nu}(1)=0.
\end{equation}

Finally, we can calculate the sum in the left-hand side
in~\eqref{eqAtheta5}. To do so, we note the following:
\begin{enumerate}
\item[(a)]  in~\eqref{eqAtheta11}, there are only terms with
$\nu^2,\nu^4,\dots,\nu^{2j}$ (since $D_{km}=0$ for $k+m$ odd);

\item[(b)] in~\eqref{eqAtheta11}, the terms with negative powers
of $(1-y)$   cancel because $\cT_\nu^{(j)}(y)$ is a polynomial
(hence, contains only nonnegative powers of $(1-y)$);

\item[(c)] after summation with respect to $\nu$,   the terms with
$\nu^4,\dots,\nu^{2j}$ will  vanish, while $\sum\limits_{\nu=1}^N
a_\nu\nu^2$ will yield 1, due to~\eqref{eqanu};

\item[(d)] the terms with positive powers of $(1-y)$ and
$\cU_\nu(y)$ will vanish at $y=1$ due to~\eqref{eqAtheta12}.
\end{enumerate}

Thus, only the terms corresponding to $(k,m)=(1,1)$ and
$(k,m)=(2,0)$ remain:
\begin{equation}\label{eqAtheta13}
(-1)^{j-1}\sum\limits_{\nu=1}^N a_\nu
\cT_\nu^{(j)}(1)=(-1)^{j-1}(-1)^j\left( D_{11}\sum\limits_{n,l\ge
0,\ n+l=j-1} C_{1l} E_{1n}+D_{20}\sum\limits_{n,l\ge 0,\ n+l=j-1}
C_{2l}E_{0n}\right).
\end{equation}
Since $u(\theta)=\cos\theta$, it follows from the first equality
in~\eqref{eqAtheta8} that $D_{11}$ and $D_{20}$ are the leading
order coefficients in the Taylor expansions of $-\sin\theta$ and
$-\cos\theta$, respectively. Hence, $D_{11}=D_{20}=-1$
and~\eqref{eqAtheta13} yields
$$
(-1)^{j-1}\sum\limits_{\nu=1}^N a_\nu \cT_\nu^{(j)}(1)=
 \sum\limits_{n,l\ge 0,\ n+l=j-1} C_{1l}
E_{1n}+ \sum\limits_{n,l\ge 0,\ n+l=j-1} C_{2l}E_{0n}>0,
$$
which completes the proof of\footnote{The left-hand side
in~\eqref{eqAtheta5} can be also found, using the following
observation. In the proof, we have shown that,
in~\eqref{eqAtheta6'}, only the terms with $u'$ and $u''$ are
relevant. Furthermore, in the expansion of $u'(\nu g(y))$ only the
term $-\nu g(y)$ is relevant, and in the expansion of $u''(\nu
g(y))$ only the term $-1$ is relevant. Therefore, we would obtain
the same result if we replaced $u(\theta)$ by the function
$-\theta^2/2$ and deleted all the terms with negative powers of
$(1-y)$ after the respective expansions in the end. Therefore,
$$
(-1)^{j-1}\sum\limits_{\nu=1}^N a_\nu \cT_\nu^{(j)}(1)=
(-1)^j\sum\limits_{\nu=1}^N a_\nu
\dfrac{d^j}{dy^j}\left(\dfrac{\nu^2 \arccos^2
y}{2}\right)\Big|_{y=1}=\dfrac{(-1)^j}{2}
\dfrac{d^j}{dy^j}(\arccos^2 y)\Big|_{y=1}
$$
because $\dfrac{d^j}{dy^j}(\arccos^2 y)=2v(y)v'(y)$ already has no
negative powers of $(1-y)$ due to~\eqref{eqAtheta6}.
Interestingly,  the left-hand side in~\eqref{eqAtheta5} does not
depend on $N$.}~\eqref{eqAtheta5} and
  assertion 1 in the lemma.

  Assertion 2 follows from assertion 1.

  The Taylor expansion in assertion 3 follows from the Taylor expansions of
$\cos\nu\theta$ and relations~\eqref{eqanu}.

\section{One identity}\label{appendGamma0}

In this section, we prove a result, which we need for
Theorems~\ref{thAsympfGeneral} and~\ref{thAsympgGeneral} in the
case $N=1$.

\begin{lemma}\label{lGammaIs0}
Let
\begin{equation*}
\Omega(x):=\dfrac{1}{2\pi}\left(\,\int\limits_{-\pi}^\pi
\left(\dfrac{\cos(x\xi)}{2(1-\cos\xi)}-\dfrac{1}{\xi^2}\right)d\xi
- \dfrac{2}{\pi} \right) + \dfrac{x}{2},\quad x\in\bbZ,\ x\ge 0.
\end{equation*}
 Then $\Omega(x)=0$ for all integer $x\ge 0$.
\end{lemma}
\proof {\bf Step 1.} First, we show that
\begin{equation}\label{eqlGammaIs0_1}
\Omega(x)-\Omega(x+1)=0\quad\text{for all integer } x\ge 0.
\end{equation}
 We have
$$
\Omega(x)-\Omega(x+1)= I(x)  - \dfrac{1}{2},
$$
where
$$
I(x):=\dfrac{1}{4\pi}\int\limits_{-\pi}^{\pi}\dfrac{\cos(x\xi)-\cos((x+1)\xi)}{
1-\cos\xi }\,d\xi.
$$

If $x=0$, then, obviously, $I(x)=1/2$.

Assume that $x\ge 1$. Using the formula
$$
\cos(x+1)\xi=\cos x\xi \cos\xi-\sin x\xi\sin\xi,
$$
we obtain
\begin{equation}\label{eqLsinntheta}
\begin{aligned}
I(x) =\dfrac{1}{4\pi}\int\limits_{-\pi}^{\pi}\dfrac{\sin
x\xi\sin\xi}{ 1-\cos\xi }\,d\xi +
\dfrac{1}{4\pi}\int\limits_{-\pi}^{\pi}\cos
x\xi\,d\xi=\dfrac{1}{4\pi}\int\limits_{-\pi}^{\pi}\dfrac{\sin
x\xi\sin\xi}{ 1-\cos\xi }\,d\xi.
\end{aligned}
\end{equation}

For $x=1$, we have $I(x)=1/2$. To prove~\eqref{eqlGammaIs0_1}, it
remains to show that the right-hand side of~\eqref{eqLsinntheta}
does not depend on $x\ge 1$. Using the formula $ \sin(x+1)\xi=\sin
x\xi \cos\xi+\cos x\xi\sin\xi, $ we obtain
$$
\begin{aligned}
\dfrac{(\sin (x+1)\xi-\sin x\xi)\sin\xi}{1-\cos\xi} &=-\sin
x\xi\sin\xi+\dfrac{\cos x\xi\sin^2\xi}{1-\cos\xi}\\
&=-\sin x\xi\sin\xi+\cos x\xi+\cos x\xi \cos\xi.
\end{aligned}
$$
For $x=1$, the integral of the right-hand side vanishes due to
direct calculation, while for $x\ge 2$, it vanishes due to the
orthogonality of the systems $\{\sin x\xi\}_{x\in\bbZ}$ and
$\{\cos x\xi\}_{x\in\bbZ}$, respectively, in $L_2(-\pi,\pi)$.

{\bf Step 2.} It remains to show that $\Omega(0)=0$. We have
\begin{equation}\label{eqlGammaIs0_2}
\Omega(0)= \dfrac{1}{\pi}\lim\limits_{\sigma\to
0}\left(I_1(\sigma)+I_2(\sigma)-\dfrac{1}{\pi}\right),
\end{equation}
where
\begin{align}
 I_1(\sigma)&:=\int\limits_{\sigma}^\pi
\dfrac{d\xi}{2(1-\cos\xi)}=\dfrac{1}{4}\int\limits_{\sigma}^{\pi}
\dfrac{d\xi}{\sin^2 (\xi/2)}=\dfrac{1}{2}\cot(\sigma/2)=\dfrac{1}{\sigma}+O(\sigma)  , \label{eqlGammaIs0_3} \\
I_2(\sigma)&:=-\int\limits_{\sigma}^\pi
\dfrac{d\xi}{\xi^2}=\dfrac{1}{\pi}-\dfrac{1}{\sigma}.\label{eqlGammaIs0_4}
\end{align}
Combining~\eqref{eqlGammaIs0_2}, \eqref{eqlGammaIs0_3},
and~\eqref{eqlGammaIs0_4} yields $\Omega(0)=0$.
\endproof

{\bf Acknowledgement.} The author expresses his gratitude to
Sergey Tikhomirov for numerous discussions. The research was
supported by DFG through the Heisenberg fellowship and the Collaborative Research Center 910: ``Control of
self-organizing nonlinear systems'' (Germany).

\end{document}